\documentclass[11pt]{article}

\input xy

\xyoption{all}

\usepackage{amssymb,amsbsy,amsmath,graphicx,epsfig,times}
\usepackage{hyperref}

\newtheorem{thm}{Theorem}[section]
\newtheorem{lem}[thm]{Lemma}
\newtheorem{prop}[thm]{Proposition}
\newtheorem{cor}[thm]{Corollary}

\newtheorem{dfn}[thm]{Definition}

\newcommand{\bs}[1]{\boldsymbol{#1}}
\renewcommand{\bf}[1]{\mathbf{#1}}
\renewcommand{\rm}[1]{\mathrm{#1}}
\renewcommand{\cal}[1]{\mathcal{#1}}

\newcommand{\bbC}{\mathbb{C}}

\newcommand{\bbN}{\mathbb{N}}

\newcommand{\bbQ}{\mathbb{Q}}
\newcommand{\bbR}{\mathbb{R}}
\newcommand{\bbT}{\mathbb{T}}
\newcommand{\bbZ}{\mathbb{Z}}

\newcommand{\bfF}{\mathbf{F}}

\newcommand{\B}{\mathcal{B}}
\newcommand{\C}{\mathcal{C}}

\newcommand{\F}{\mathcal{F}}

\renewcommand{\P}{\mathcal{P}}

\newcommand{\frH}{\mathfrak{H}}

\newcommand{\frh}{\mathfrak{h}}

\newcommand{\G}{\Gamma}
\renewcommand{\O}{\Omega}

\renewcommand{\L}{\Lambda}

\renewcommand{\a}{\alpha}

\renewcommand{\l}{\lambda}
\renewcommand{\o}{\omega}

\newcommand{\id}{\mathrm{id}}
\newcommand{\Cay}{\mathrm{Cay}}

\renewcommand{\hat}[1]{\widehat{#1}}

\newcommand{\into}{\hookrightarrow}

\newcommand{\fin}{\nolinebreak\hspace{\stretch{1}}$\lhd$}
\newcommand{\qed}{\nolinebreak\hspace{\stretch{1}}$\Box$}

\newcommand{\actson}{\curvearrowright}

\begin{document}

\title{Rational group ring elements with kernels having irrational dimension}
\author{Tim Austin\\ \\ \small{Courant Institute of Mathematical Sciences,}\\ \small{New York University,}\\ \small{NY
10012, USA}\\ \small{\texttt{tim@cims.nyu.edu}}}

\date{}

\maketitle

\begin{abstract}
We prove that there are examples of finitely generated groups $\G$
together with group ring elements $Q \in \bbQ\G$ for which the von
Neumann dimension $\dim_{L\G}\ker Q$ is irrational, so (in
conjunction with other known results) answering a question of
Atiyah.
\end{abstract}


\parskip 0pt

\tableofcontents

\parskip 7pt

\section{Introduction}

Given a countable discrete group $\G$, we write $\bbQ\G$ and
$\bbC\G$ respectively for its rational and complex group rings,
$\l:\G \actson \ell^2(\G)$ for the Hilbertian completion of its left
regular representation and $L\G$ for the resulting group von Neumann
algebra, which may be obtained by completing $\l(\bbC\G)$ in the
weak operator topology of $\B(\ell^2(\G))$. Henceforth we will
generally identify $\bbQ\G$ and $\bbC\G$ with their images in $L\G$
under $\l$. In this setting we can define the von Neumann dimension
of any closed $L\G$-submodule of $\ell^2(\G)$; we assume familiarity
with this notion, referring the reader to the book of
L\"uck~\cite{Luc02} for an introduction. We will address the
following classical question:

\begin{quote}
Do there exist $\G$ and $Q \in \bbQ\G$ for which the von Neumann
dimension $\dim_{L\G}\ker Q$ is irrational?
\end{quote}

This is known to be equivalent to the problem posed by Atiyah of
constructing a cocompact free proper $\G$-manifold without boundary
that has irrational $L^2$-Betti numbers (originally formulated as
problem (iii) on page 72 of~\cite{Ati76}). This equivalence is
proved in Lemma 10.5 of L\"uck~\cite{Luc02}: in particular, it is
proved that given any $\G$ and $A \in \bbQ\G$ one can construct a
cocompact free proper $\G$-manifold one of whose $L^2$-Betti numbers
is equal to $\dim_{L\G}\ker A$. We will henceforth restrict our
attention to the purely group-theoretic version of the problem.  A
much more thorough discussion of this question is contained in
L\"uck's~\cite{Luc02} Chapter 10, and a discussion of its relation
to questions of computability can be found in section 8.A$_4$ of
Gromov's essay in~\cite{NibRol93}.

A stronger version of the question, asking whether in fact
$\dim_{L\G}\ker Q$ must always lie in the additive subgroup
$\rm{fin}^{-1}(\G) \leq \bbQ$ generated by the inverses of the
orders of the finite subgroups of $\G$, is now known to be false
from the work~\cite{GriZuk01} of Grigorchuk and \.Zuk (see also the
article~\cite{GriLinSchZuk00} of Grigorchuk, Linnel, Schick and
\.Zuk), who have shown that the lamplighter group $\bbZ_2\wr \bbZ$
is a counterexample: all of its finite subgroups have order that is
a member of $2^\bbZ$, but a natural finitely-supported operator with
integer coefficients on the group (in fact, a rational multiple of a
Markov operator) has an eigenspace with von Neumann dimension
$\frac{1}{3}$.

A new and quite elementary treatment of this fact has now been given
by Dicks and Schick in~\cite{DicSch02}, and in this work we will
adapt some of their calculations to provide a family of examples
answering the original question about irrational dimensions, as
formulated above. In order to state our main theorem precisely we
first need a little notation.

We write $\bfF_n$ to denote the free group on $n$ generators, $s_1$,
$s_2$, \ldots, $s_n$ for those generators themselves, $S =
\{s_1^{\pm 1},s_2^{\pm 1},\ldots,s_n^{\pm 1}\}$ for the
corresponding symmetric generating set and $e$ for the identity
element of $\bfF_n$. To these data are associated the Cayley graph
$\Cay(\bfF_n,S)$ with vertex set $\bfF_n$ and edge set $\{\{g,gs\}:\
g\in \bfF_n,\,s \in S\}$, which is simply a $2n$-regular infinite
tree. Here and later in the paper we will use mostly standard
graph-theoretic terminology in relation to $\Cay(\bfF_n,S)$, as
described, for instance, in Chapter I of Bollob\'as~\cite{Bol98}.
Given a subset $A \subset \bfF_n$ we will write $\Cay(\bfF_n,S)|_A$
for the induced subgraph of $\Cay(\bfF_n,S)$ on the set of vertices
$A$, and
\[\partial A := A\cdot S\setminus A\]
for the \textbf{boundary} of $A$ in $\Cay(\bfF_n,S)$.  A
\textbf{path} in $\Cay(\bfF_2,S)$ is a subset $P =
\{g_0,g_1,\ldots,g_\ell\} \subset \bfF_2$ with $g_{i+1} \in g_iS$
for every $i\leq \ell - 1$ and with all the $g_i$s distinct, and in
this case the \textbf{length} of the path is $\ell$. We denote by
$\rho$ the left-invariant word metric on $\bfF_2$, which is simply
the graph distance arising from $\Cay(\bfF_n,S)$, and will sometimes
refer to $\rho(e,g)$ as the \textbf{length} of an element $g \in
\bfF_n$.  Given $g \in \bfF_n$ and $r \geq 0$ we let $B(g,r) := \{h
\in \bfF_n:\ \rho(g,h) \leq r\}$ be the closed ball of radius $r$
around $g$ in $\Cay(\bfF_n,S)$, and more generally given $A\subseteq
\bfF_n$ we let $B(A,r) := \bigcup_{g\in A}B(g,r)$ be its
\textbf{radius-$r$ neighbourhood}.

In addition we write $\bbZ_2$ to denote the cyclic group of order
$2$, and $\bbZ_2^{\oplus I}$ (respectively $\bbZ_2^I$) to denote the
direct sum (respectively direct product) of a family of copies of
$\bbZ_2$ indexed by some other set $I$. We will usually denote
members of $\bbZ_2^{\oplus I}$ by lowercase bold letters such as
$\bf{w} = (w_i)_{i\in I}$, and will write $\delta_i$ for the
distinguished element of $\bbZ_2^{\oplus I}$ that takes the value $1
\in \bbZ_2$ at $i$ and $0$ elsewhere.

The main result of this paper is the following.

\begin{thm}\label{thm:main}
Let the space $\P(\bbN)$ of subsets of $\bbN$ be endowed with the
lexicographic ordering. There are parameterizations
\[\P(\bbN) \ni I \mapsto V_I\leq \bbZ_2^{\oplus
\bfF_2}\] of a family of subgroups that are invariant under the
left-coordinate-translation action of $\bfF_2$ and
\[I \mapsto Q_I \in \bbQ((\bbZ_2^{\oplus \bfF_2}/V_I) \rtimes \bfF_2)\]
of a family of rational group ring elements such that the associated
map
\[\P(\bbN)\to \bbR:I\mapsto \dim_{L((\bbZ_2^{\oplus \bfF_2}/V_I)\rtimes
\bfF_2)}\ker(Q_I - 4)\] is strictly increasing, where $\bfF_2\actson
\bbZ_2^{\oplus \bfF_2}/V_I$ by left-coordinate-translation.
\end{thm}

Since a strictly increasing map is an injection, the image in $\bbR$
of $\P(\bbN)$ under this map must be uncountable, and so we
immediately obtain the following.

\begin{cor}
For some left-translation-invariant subspace $V \leq \bbZ_2^{\oplus
\bfF_2}$, the finitely-generated group $(\bbZ_2^{\oplus
\bfF_2}/V)\rtimes \bfF_2$ admits a group ring element with rational
coefficients whose kernel has irrational (and even transcendental)
von Neumann dimension. \qed
\end{cor}

The main innovation of this paper is to exploit the freedom in the
choice of the subgroup $V$ above in order to obtain a large family
of von Neumann dimensions, some of which must then be irrational,
rather than trying to find one single example of a group and group
ring element and compute the von Neumann dimension of its kernel
explicitly. It is this idea that we will make precise in obtaining
the family of examples promised in Theorem~\ref{thm:main}.  A
similar instance of exploiting this freedom in the choice of $V$ to
produce an example of a group with interesting properties appeared
recently in~\cite{Aus--amenablepoorcompression}, and the present
paper was indirectly motivated by that one.

\subsubsection*{Remark}

Since a version of the present paper first appeared online~\cite{Aus--Atiyahexample}, works of Pichot, Schick and \.{Z}uk~\cite{PicSchZuk10} and  Grabowski~\cite{Gra10} have used a similar underlying construction to produce a range of related examples.  By incorporating several non-trivial new ideas, those examples can be made simpler and more explicit than in the present paper, and can be arranged to have various additional properties such as amenability and finite presentation.  Even more concrete examples have been given by Lehner and Wagner in~\cite{LehWag10}, using an extension of the ideas from~\cite{DicSch02}.

\subsubsection*{Acknowledgements}

I am grateful to Dimitri Shlyakhtenko for explaining to me much of
the background to Atiyah's question, and to Wolfgang L\"uck, Peter
Linnell and Thomas Schick for insightful suggestions on improving
this paper.

\section{Some preliminary manipulations}\label{sec:prelims}

In this section we let $\L$ be any discrete group and $U$ any
discrete Abelian group equipped with a left action $\a:\L \actson U$
by automorphisms (so $\a^{gh} = \a^g\circ \a^h$).  From these we form the
semidirect product $U\rtimes_\a \L$ as the set-theoretic Cartesian
product $U\times \L$ with the multiplication
\[(u,g)\cdot (w,h) := (\a^{h^{-1}}(u) + w,gh).\]
Let $\hat{U}$ be the compact dual group of $U$, and let $m_{\hat{U}}$ be its Haar probability measure.

We now describe an identification of the left regular action
\[\l:\big(\bbC(U\rtimes_\a \L) \subset L(U\rtimes_\a \L)\big)\actson \ell^2(U\rtimes_\a \L)\]
that will prove convenient later.

The point is simply that the Fourier transform sets up a unitary isomorphism
\[\F:\ell^2(U) \stackrel{\cong}{\longrightarrow} L^2(m_{\hat{U}})\]
\[\F f(\chi) := \sum_{u \in U}\langle u,\chi\rangle f(u),\]
and now since $U \rtimes_\a \L$ is set-theoretically simply equal to
$U \times \L$, we have also
\[\F\otimes \rm{Id}_{\ell^2(\L)}:\ell^2(U\rtimes_\a\L) \cong \ell^2(U)\otimes \ell^2(\L) \stackrel{\cong}{\longrightarrow} L^2(m_{\hat{U}})\otimes \ell^2(\L) \cong L^2(m_{\hat{U}}\otimes \#_\L),\]
where we write $\#_S$ to denote the counting measure on a set $S$.
Let $\hat{\a}:\L\actson \hat{U}$ be the Pontrjagin adjoint action of
$\a$ defined by the relation
\[\langle u,\hat{\a}^g(\chi)\rangle := \langle \a^{g^{-1}}(u),\chi\rangle,\]
and recall that the duality $\langle\cdot,\cdot\rangle:U\times \hat{U}\to \bbT$ establishes the Pontrjagin isomorphism $U \cong \hat{\hat{U}}$.

As is standard, the isomorphism $\F$ of Hilbert spaces now defines an isomorphism of actions
\[\Bigg(\l:\begin{array}{c}\bbC(U\rtimes_\a\L)\\ \cap\\ L(U\rtimes_\a\L)\end{array}\actson \ell^2(U\rtimes_\a\L)\Bigg) \stackrel{\cong}{\longrightarrow} \Bigg(\pi:\begin{array}{c}\bbC(U\rtimes_\a\L)\\ \cap\\ L^\infty(m_{\hat{U}})\rtimes_{\hat{\a}}\L\end{array}\actson L^2(m_{\hat{U}}\otimes\#_\L)\Bigg),\]
where $\pi$ is the left regular action of the group measure space
von Neumann algebra $L^\infty(m_{\hat{U}})\rtimes_{\hat{\a}} \L$ on
$L^2(m_{\hat{U}}\otimes \#_\L)$, and within
$\cal{U}\big(L^\infty(m_{\hat{U}})\rtimes_{\hat{\a}} \L\big)$ we
identify copies of $U \cong \hat{\hat{U}}$ and $\L$ that together
generate a copy of $U\rtimes_\a \L$ acting by
\[\pi(u,e)f(\chi,g) = \langle \a^{g^{-1}}(u),\chi\rangle f(\chi,g) = (M_{\langle u,\cdot\rangle}f)(\chi,g)\quad\quad \hbox{for}\ (\chi,g)\in \hat{U}\times \L,\]
where $M_F$ denotes twisted pointwise multiplication
\[M_F f(\chi,g) := F(\hat{\a}^g(\chi))f(\chi,g)\]
by a function $F \in L^\infty(m_{\hat{U}})$, and
\[\pi(0,h)f(\chi,g) = f(\chi,h^{-1}g) =: T^hf(\chi,g)\quad\quad \hbox{for}\ (\chi,g) \in \hat{U}\times \L,\]
so this is still a translation operator.  If $W \subseteq \hat{U}$
is a Borel subset we will generally write $M_W$ in place of
$M_{1_W}$.

We may double-check that the above specifications do combine to give
an action of $U\rtimes_\a\L$ through the following commutation
relation:
\begin{eqnarray*}
(T^{h^{-1}}\circ M_F\circ T^h)f(\chi,g) &=& ((M_F\circ
T^h)f)(\chi,hg)\\
&=& F(\hat{\a}^{hg}(\chi))\cdot (T^hf)(\chi,hg)\\
&=& (F\circ\hat{\a}^h)(\hat{\a}^g(\chi))\cdot f(\chi,g)\\
&=& M_{F\circ \hat{\a}^h}f(\chi,g)
\end{eqnarray*}
for $F \in L^\infty(m_{\hat{U}})$ and $h \in \L$.

These manipulations lead to a simple identification between group
von Neumann algebras $L(U\rtimes_\a \L)$ and group measure space
algebras $L^\infty(m_{\hat{U}})\rtimes_{\hat{\a}}\L$ corresponding
to dynamical systems $\hat{\a}:\L\actson \hat{U}$ of algebraic
origin.  In the case $\L = \bbZ^d$ for $d \geq 2$ such dynamical
systems are known to exhibit a wide variety of interesting behaviour
(see, in particular, the monograph~\cite{Sch95} of Schmidt), and in
recent years the analysis of such systems for certain non-Abelian
$\L$ has also begun to make headway (see, for instance, the
paper~\cite{DenSch07} of Deninger and Schmidt and the further
references given there). In the present paper we make our own modest
appeal to this dynamical picture of semidirect products with Abelian
kernel, and it would be interesting to explore whether insights from
that field could be used to drive other constructions in geometric
group theory in the future.

The above shows that to study the von Neumann algebra properties of
$\l(\bbQ(U\rtimes_\a \L))$ (turning our attention now to the
rational group ring) we may equivalently consider
$\pi(\bbQ(U\rtimes_\a \L))$, whose members may all be put into the
form
\[\sum_{i=1}^n T^{g_i}\circ M_{\phi_i}\]
with $g_i \in \L$ for each $1 \leq i \leq n$ and each $\phi_i \in
C(\hat{U})$ being a trigonometric polynomial with rational
coefficients (that is, a finite $\bbQ$-linear combination of
characters) on $\hat{U}$. It is this form for our operators that
will be most convenient for the proof of Theorem~\ref{thm:main}.

We will henceforth apply the above manipulations in case $\L=
\bfF_2$, and will specialize to groups $U$ of the form
$\bbZ_2^{\oplus \bfF_2}/V$ for some left-translation-invariant
subgroup $V \leq \bbZ_2^{\oplus \bfF_2}$, equipped with the left
translation action
\[\a^g((u_h)_{h \in \bfF_2} + V) = (u_{g^{-1}h})_{h \in \bfF_2} + V.\]

In this case the Pontrjagin duals obey the relations
\[\widehat{\bbZ_2^{\oplus \bfF_2}} \cong \hat{\bbZ_2}^{\bfF_2}\cong \bbZ_2^{\bfF_2}\]
and
\[\widehat{\bbZ_2^{\oplus \bfF_2}/V} \cong V^\perp := \big\{\bs{\chi} \in \bbZ_2^{\bfF_2}:\ \langle \bf{v},\bs{\chi}\rangle = 0\ \forall \bf{v} \in V\big\}.\]

We recognize $\hat{\a}:\bfF_2 \actson V^\perp$ as the subshift of
the left-acting topological Bernoulli shift $\bfF_2\actson
\bbZ_2^{\bfF_2}$ defined by the relations of annihilating all
members of $V$.  Let us also note for future reference that with
these conventions, if $\bs{\chi} = (\chi_h)_{h\in \bfF_2} \in
\bbZ_2^{\bfF_2}$ and $g \in \bfF_2$, then regarding $\bs{\chi}$ as a
colouring of $\rm{Cay}(\bfF_2,S)$ by elements of $\bbZ_2$, the point
$\hat{\a}^{g^{-1}}(\bs{\chi})$ is obtained by shifting that
colouring by the graph automorphism of $\rm{Cay}(\bfF_2,S)$ that
moves the point $g$ to the origin and respects the directions of all
the edges.

Note also that in this case the rational trigonometric polynomials
on $\hat{U}$ are easily seen to be those functions on $V^\perp$ that
are restrictions of functions on $\bbZ_2^{\bfF_2}$ that depend on
only finitely many coordinates and that take only rational values
(using the fact that characters on groups of the form $\bbZ_2^I$
take only the values $\pm 1$, so in particular are all
rational-valued), and so henceforth we will freely work with such
functions when specifying members of $\pi(\bbQ((\bbZ_2^{\oplus
\bfF_2}/V)\rtimes_\a\bfF_2))$ of interest. We will also now work
only with left-translation actions of $\bfF_2$ such as the above,
and so will usually omit their explicit mention from our notation.

\section{Introduction of the operators}\label{sec:intro-operators}

\subsection{Construction}

We now introduce certain members of
$\pi(\bbQ((\bbZ_2^{\oplus \bfF_2}/V)\rtimes\bfF_2))$.  These will
take the form
\[Q = \sum_{s \in S}T^{s^{-1}}\circ (M_{F_s} + M_{G_s\circ \hat{\a}^{s^{-1}}})\]
where $F_s,G_s:\bbZ_2^{\bfF_2} \to \bbQ$ for $s \in S$ depend only
on some finite patch of coordinates around $e \in \bfF_2$.  Note
that in considering the above operator as a member of
$\pi(\bbQ((\bbZ_2^{\oplus \bfF_2}/V)\rtimes\bfF_2))$, we are
implicitly regarding the above as a shorthand for
\[\sum_{s \in S}T^{s^{-1}}\circ (M_{F_s|_{V^\perp}} + M_{G_s|_{V^\perp}\circ \hat{\a}^{s^{-1}}});\]
we will generally overlook this notational detail in the following.

The rather redundant form in which $Q$ has been written above, with
a sum of two terms of the form $M_F$ for each $s\in S$, is
convenient in view of the following simple calculation.

\begin{lem}\label{lem:self-adjoint} If $F_s = G_{s^{-1}}$ for every $s\in S$ then $Q$ is
self-adjoint.
\end{lem}

\textbf{Proof}\quad Since $M_F$ is self-adjoint whenever $F$ takes
real values and $(T^s)^\ast = T^{s^{-1}}$, we deduce from the
commutator relation for these operators, the symmetry of $S$ and our
assumption that
\begin{eqnarray*}
Q^\ast &=& \sum_{s \in S}(M_{F_s} + M_{G_s\circ
\hat{\a}^{s^{-1}}})\circ T^s\\ &=& \sum_{s \in S}T^s\circ (M_{F_s
\circ \hat{\a}^s} + M_{G_s \circ \hat{\a}^{s^{-1}}\circ
\hat{\a}^s})\\ &=& \sum_{s \in S}T^s\circ (M_{G_{s^{-1}} \circ
\hat{\a}^s} + M_{F_{s^{-1}}}) = Q.
\end{eqnarray*}
\qed

Most of this section will be concerned with the choice of $F_s$ and
$G_s$, which will be pivotal for what follows. We will choose
functions that depend only on coordinates in the ball $B(e,100)$.
Heuristically, the values of $F_s(\bs{\chi})$ will depend on
different features of the level-set $\bs{\chi}^{-1}\{0\}$ describing
in what ways it locally resembles a path in $\Cay(\bfF_2,S)$, what
that path looks like, and whether it contains $e$. To explain this
we first make the following useful definitions.

\begin{dfn}[Small horizontal doglegs]
A (finite or infinite) path $P \subset \Cay(\bfF_2,S)$ contains a
\textbf{small horizontal dogleg} if it contains a subset of the form
\begin{multline*}
\{gs_2^{\eta'},g,gs_1^\eta,gs_1^{2\eta},\ldots,gs_1^{\ell\eta},gs_1^{\ell\eta}
s_2^{\eta''} \}\quad\quad\hbox{for some}\ g \in \bfF_2,\, \ell \in
\{1,2,\ldots,9\},\\ \eta',\eta,\eta'' \in \{-1,1\}
\end{multline*}
or of the form
\begin{multline*}
\{g,gs_1^\eta,gs_1^{2\eta},\ldots,gs_1^{\ell\eta},gs_1^{\ell\eta}
s_2^{\eta''}\}\quad\quad\hbox{for some}\ g \in \bfF_2,\, \ell \in
\{1,2,\ldots,9\},\\ \eta,\eta''\in \{-1,1\},\, \hbox{with}\ g\
\hbox{an end-point of}\ P.
\end{multline*}
(note that only the first of these cases really fits the term
`dogleg'). Otherwise $P$ contains \textbf{no small horizontal
doglegs}.  In either case we refer to the further subset
$\{g,gs_1^\eta,gs_2^{2\eta},\ldots,gs_1^{\ell\eta}\}$ as the
\textbf{main segment} of the dogleg.
\end{dfn}

\begin{dfn}[Locally good points]
A point $\bs{\chi} \in \bbZ_2^{\bfF_2}$ is \textbf{locally good} if
\begin{enumerate}
\item $\bs{\chi}^{-1}\{0\}\cap B(e,10)$ is a path in
$\Cay(\bfF_2,S)|_{B(e,10)}$ that contains $e$ and has length at
least $10$ (that is, it connects $e$ with some point of $\partial
B(e,9) \subset B(e,10)$),
\item there is no small horizontal dogleg in the path
$\bs{\chi}^{-1}\{0\}\cap B(e,10)$ whose main segment lies within
$B(e,9)$, and
\item for every $g \in \bs{\chi}^{-1}\{0\}\cap B(e,10)$ we also have
that $\bs{\chi}^{-1}\{0\}\cap B(g,10)$ is a path in
$\Cay(\bfF_2,S)|_{B(g,10)}$ containing no small horizontal doglegs
with main segment contained in $B(g,9)$.
\end{enumerate}
\end{dfn}

The second part of the above definition is very important.  It
places rather severe restrictions on which paths can appear as
$\bs{\chi}^{-1}\{0\}\cap B(e,10)$ if $\bs{\chi}$ is locally good:
insofar as a path in $\Cay(\bfF_2,S)$ is made up of a concatenation
of `horizontal' segments (with steps given by $s_1^{\pm 1}$) and
`vertical' segments (with steps given by $s_2^{\pm 1}$), this
condition tells us that while the maximal vertical segments that
appear in $\bs{\chi}^{-1}\{0\}\cap B(e,10)$ may be of any length,
this path may not contain any maximal horizontal segments that lie
properly inside $B(e,10)$ and have length less than $10$.  It
follows that if a maximal horizontal segment lies properly inside
$B(e,10)$ (that is, the end-points of that segment also visibly lie
inside $B(e,10)$), then it must contain $e$ as an interior point and
extend to points $s_1^a$ and $s_1^{-b}$ with $a,b\geq 1$ and $a +
b\geq 10$, before being permitted to make at most one more
horizontal-vertical-horizontal `dogleg' before leaving $B(e,10)$ on
either side of $e$.  Moreover, the last condition of the above
definition ensures that not only does the vertex $e$ see this highly
constrained behaviour in its radius-$10$ neighbourhood, but also all
of its neighbours inside this path and at distance at most $10$ see
this behaviour in their radius-$10$ neighbourhoods. This rather
peculiar restriction on the kinds of path we allow will be pivotal
at exactly one point below (Corollary~\ref{cor:always-get-zero}),
where it will restrict a certain sum over paths to terms that
possess some additional helpful properties.

We will give a definition of $F_s$ (and then set $G_s = F_{s^{-1}}$)
that uses the above notion, but we first define another auxiliary
function $F^\circ_s$.

\begin{dfn}\label{dfn:FScirc}
The function $F_s^\circ:\bbZ_2^{\bfF_2}\to \bbQ$ is defined
according to the following four cases:
\begin{itemize}
\item $F^\circ_s(\bs{\chi}) := 1$ if $\bs{\chi}$ is locally good
and $e$ and $s$ are both interior points of the path
$\bs{\chi}^{-1}\{0\}\cap B(e,10)$;
\item $F^\circ_s(\bs{\chi}) := 2$ if $\bs{\chi}$ is locally good, $e$ is
an interior point of the path $\bs{\chi}^{-1}\{0\}\cap B(e,10)$ and
$s$ is its end-point; or if $\bs{\chi}$ is locally good and the path
$\bs{\chi}^{-1}\{0\}\cap B(e,10)$ contains both $e$ and also some $t
\in S\setminus\{s,s^{-1}\}$, but does not contain $s$;
\item $F^\circ_s(\bs{\chi}) := \frac{1}{100}$ if $\bs{\chi}$ is \emph{not} locally
good, but we do have that $e \in \bs{\chi}^{-1}\{0\}$ and that the
translate $\hat{\a}^{s^{-1}}(\bs{\chi})$ is locally good.
\item $F^\circ_s(\bs{\chi}) := 0$ otherwise.
\end{itemize}
\end{dfn}

\textbf{Remarks}\quad\textbf{1.}\quad In particular,
$F^\circ_s(\bs{\chi}) = 0$ unless $e \in \bs{\chi}^{-1}\{0\}$ and
$\bs{\chi}^{-1}\{0\}\cap B(e,10)$ is a path in
$\Cay(\bfF_2,S)|_{B(e,10)}$, and given these conditions the exact
value of $F^\circ_s(\bs{\chi})$ is determined by a further
sub-classification.

\quad\textbf{2.}\quad Let us draw attention to the quirk that if
$\bs{\chi}$ is locally good but $e$ is an \emph{end-point} of the
path $\bs{\chi}^{-1}\{0\}\cap B(e,10)$ with neighbour $s$ also lying
in this path, then
\begin{eqnarray*}
F_s(\bs{\chi}) &=& F_{s^{-1}}(\bs{\chi}) = 0,\\
G_s(\hat{\a}^{s^{-1}}(\bs{\chi})) &=&
F_{s^{-1}}(\hat{\a}^{s^{-1}}(\bs{\chi})) = 2\\
\hbox{and}\quad F_t(\bs{\chi}) &=& 2\quad \hbox{for}\ t \in
S\setminus\{s,s^{-1}\}.
\end{eqnarray*}
This slightly tricky case will give rise to a useful simplification
later.

\quad\textbf{3.}\quad In the third case above we must have that
$\bs{\chi}^{-1}\{0\}\cap B(e,10)$ is a path containing $e$, so this
case can arise only because there is some point $g \in
\bs{\chi}^{-1}\{0\}$ that lies at distance $10$ from $e$ and $11$
from $s$, such that $g$ also lies at distance $10$ from some `bad'
feature of $\bs{\chi}^{-1}\{0\}$ --- a fork, a distinct connected
component, or a small horizontal dogleg visible in its entirety
--- so that some other condition in the definition of `locally good'
is violated. It is easy to see that in this scenario there can be
only one such $s$, since if $s' \in S$ were another then the path
$\bs{\chi}^{-1}\{0\}\cap B(e,10)$ would have to contain both $s$ and
$s'$, and so must connect them via $e$, but in this case we see
easily from the definition that if both
$\hat{\a}^{s^{-1}}(\bs{\chi})$ and $\hat{\a}^{(s')^{-1}}(\bs{\chi})$
are locally good then so is $\bs{\chi}$.

\quad\textbf{4.}\quad Of course, the particular value
$\frac{1}{100}$ employed in the third case above is not very
important; it has been chosen simply as a rational number that will
easily be shown to satisfy a certain modest algebraic condition that
we need later. \fin

We also let $G^\circ_{s^{-1}} := F^\circ_s$, and now note the
following consequence of this definition.

\begin{lem}\label{lem:connected-components}
For any $\bs{\chi} \in \bbZ_2^{\bfF_2}$ the set
\[E(\bs{\chi}) := \bigcup_{s\in S}\big\{g \in \bfF_2:\ F^\circ_s(\hat{\a}^{g^{-1}}(\bs{\chi}))\ \hbox{and}\ G^\circ_s(\hat{\a}^{s^{-1}g^{-1}}(\bs{\chi}))\ \hbox{not both}\ 0\big\}\]
is a union of connected components in $\Cay(\bfF_2,S)$ each of which
takes the form $B(P,1)\setminus T$ for some path $P$ with no small
horizontal doglegs and some set $T$ of at most two boundary points
of end-points of $P$, and any two of these connected components are
separated by a distance of at least $9$.
\end{lem}

\textbf{Proof}\quad If $F^\circ_s(\hat{\a}^{g^{-1}}(\bs{\chi}))$ and
$G^\circ_s(\hat{\a}^{s^{-1}g^{-1}}(\bs{\chi}))$ are not both zero,
then from the first remark above it follows that either $g$ is
itself a member of $\bs{\chi}^{-1}\{0\}\cap B(g,10)$ and
$\hat{\a}^{g^{-1}}(\bs{\chi})$ is locally good, so this set takes
the form of a path with no small horizontal doglegs in
$\Cay(\bfF_2,S)|_{B(g,10)}$, or $g$ is adjacent to such a point.
Clearly these paths in balls of radius $10$ patch together to form,
together with their immediate neighbourhoods, the connected
components of the given set, so each of these must be a path in
$\Cay(\bfF_2,S)$ together with all but possibly two members of its
neighbourhood (these being precisely the points such as $s^{-1}$ in
the situation described in Remark 2 above). From the definition of
$F_s$ it follows that any two distinct such paths must lie at a
distance of at least $11$ from each other (and so their radius-$1$
neighbourhoods must lie at distance at least $9$), in order that
points internal to these paths should not see foreign connected
components within their radius-$10$ neighbourhoods. \qed

\begin{cor}
The set
\begin{multline*}
W := \big\{\bs{\chi}\in \bbZ_2^{\bfF_2}:\ E(\bs{\chi})\ni e\
\hbox{but the central path of the component that}\\ \hbox{contains}\
e\ \hbox{has length}\ \leq 4\big\}
\end{multline*}
depends only on coordinates in $B(e,100)$. \qed
\end{cor}

Finally, we set
\[F_s := F^\circ_s\cdot 1_W\quad\quad G_s := F_{s^{-1}} = G^\circ_s\cdot 1_W\]
and consider the resulting operator $Q$, which by
Lemma~\ref{lem:self-adjoint} is self-adjoint.

\subsection{Decomposition into invariant subspaces}

In describing the further consequences of our choice of $F_s$ the
following terminology will prove convenient.

\begin{dfn}[Good and bad neighbourhoods]
For a given point $\bs{\chi}\in \bbZ_2^{\bfF_2}$, a ball $B(g,10)
\subset \bfF_2$ is a \textbf{good neighbourhood for $\bs{\chi}$} if
$\hat{\a}^{g^{-1}}(\bs{\chi})$ is locally good and $g$ is an
end-point of the path $\bs{\chi}^{-1}\{0\}\cap B(g,10)$.  It is a \textbf{bad neighbourhood for $\bs{\chi}$} if $\hat{\a}^{g^{-1}}(\bs{\chi})$ is not locally
good, but $g \in \bs{\chi}^{-1}\{0\}$ and for some $s \in S$ the
translate $\hat{\a}^{s^{-1}g^{-1}}(\bs{\chi})$ is locally good.
\end{dfn}

Now consider a point $\bs{\chi} \in \bbZ_2^{\bfF_2}$: either there
is some $s \in S$ such that \[F_s(\bs{\chi})\ \hbox{and}\
G_s(\hat{\a}^{s^{-1}}(\bs{\chi}))\ \hbox{are not both}\ 0,\] or
there is not. Let $C_0$ be the set of those $\bs{\chi}$ for which
there is not; this is clearly a clopen subset of $\bs{\chi}$. Our
next step will be to obtain a rather more detailed partition of the
remainder $\bbZ_2^{\bfF_2}\setminus C_0$.

Thus, suppose now that $\bs{\chi} \in \bbZ_2^{\bfF_2}\setminus C_0$,
and that $s \in S$ is such that $F_{s}(\bs{\chi})\neq 0$. It follows
that either $\bs{\chi}$ is locally good, or (if $F_{s}(\bs{\chi}) =
\frac{1}{100}$) that $e \in \bs{\chi}^{-1}\{0\}$ and
$\hat{\a}^{s^{-1}}(\bs{\chi})$ is locally good.  In either case this
requires that $\bs{\chi}^{-1}\{0\}\cap B(e,10)$ be a path with no
small horizontal doglegs that passes through $e$.

Similarly, if $G_{s}(\hat{\a}^{s^{-1}}(\bs{\chi}))\neq 0$, then
either $\bs{\chi}$ is locally good and so $\bs{\chi}^{-1}\{0\}\cap
B(e,10)$ is a path that passes through $e$, or
$\hat{\a}^{s^{-1}}(\bs{\chi})$ is locally good and
$\bs{\chi}^{-1}\{0\}\cap B(e,9)$ is a path containing $s$ but no
other member of $S\cup\{e\}$.

In either of the above cases we may pick a
unique $g_0 \in S\cup \{e\}$ that is closest to $e$ and such that
$\hat{\a}^{g_0^{-1}}(\bs{\chi})$ is locally good.

Now imagine dispatching two walkers from $g_0$ towards the two
different end-points of the path $\bs{\chi}^{-1}\{0\}\cap B(g_0,10)$
with instructions to walk in their given directions along edges that
remain in the level set $\bs{\chi}^{-1}\{0\}$ and through vertices
$g$ such that $\hat{\a}^{g^{-1}}(\bs{\chi})$ is still locally good,
until they reach either a good neighbourhood or a bad neighbourhood
for $\bs{\chi}$, where they should stop and report back to us. It
may happen that one or both of them leave $B(g_0,10)$, or that they
do not move at all.

If a walker never reaches a good or bad neighbourhood, then it
follows that the level set $\bs{\chi}^{-1}\{0\}$ that she followed
in her direction must continue to look like a path, with no
end-points, forks, small horizontal doglegs or distinct components
lying within distance $10$ of it: otherwise the walker would at some
point have stopped walking in a bad neighbourhood. Let us call this
walking-forever scenario $(\infty)$.

If the walker reaches a good neighbourhood, then she has followed a
path-like branch of $\bs{\chi}^{-1}\{0\}$ with no small horizontal
doglegs until reaching an end-point of that path, and again this
finite-length path-like branch has no other points of
$\bs{\chi}^{-1}\{0\}$ lying within distance $10$ of it. Note that
this includes the possibility that $g_0$ is an end-point of the
path, and so this walker is already in a good neighbourhood
initially. We call this ending scenario $(1)$.

The final scenario, that the walker's journey terminates in a bad
neighbourhood, may result from three different features of
$\bs{\chi}^{-1}\{0\}$: a point of this level set not connected to
the walker's path, but lying within distance $10$ of it; a fork in
the path; or a small horizontal dogleg.  In any case the walker stop
walking as soon as he reaches within distance $10$ of some further
point of his path which, in turn, can see this feature within its
radius-$10$ neighbourhood (effectively he has had a premonition of
this bad feature within distance $10$ of his own radius-$10$
horizon). This rather convoluted description is important, because
it causes this walker to stop far short of actually reaching, or
even being himself able to see, this non-path-like feature (rather
than, for example, continuing until he actually reaches a fork), and
we will find that this greatly simplifies certain enumerations
later. Note that in case $g_0 \neq e$, this includes the possibility
that this walker is dispatched from $g_0$ back towards $e$, then
reaches $e$ where this happens and stops. We call this ending
scenario $(2)$.

Finally, note also that from the definition of $F_s$ as
$F^\circ_s\cdot 1_W$, the combined distances walked by the two
walkers must be at least $5$; this rules out some annoying
degenerate scenarios, and was why we introduced the set $W$.

Now, every point $\bs{\chi} \in \bbZ_2^{\bfF_2}\setminus C_0$
results in a pair of ending scenarios, each from the set
$\{(1),(2),(\infty)\}$, according to the fates of the two walkers.
Together their route specifies some (finite or infinite) path $P
\subseteq \bs{\chi}^{-1}\{0\}$. Regarding the two walkers as
indistinguishable except by their ending scenarios, we can now
partition $\bbZ_2^{\bfF_2} \setminus C_0$ into the six (manifestly
Borel) sets $C_{a,b}$ for $a,b \in \{1,2,\infty\}$ and $a \leq b$,
where $\bs{\chi}\in C_{a,b}$ if one walker ends in scenario $(a)$
and the other in scenario $(b)$. Also, if either walker ends in a
bad neighbourhood, then we know that the path they were following
extends another $10$ steps beyond their ending position to a point
that can see bad behaviour within distance $10$ of itself, and so
including these last few steps if available defines a larger path $R
\subseteq \bs{\chi}^{-1}\{0\}$, $R \supseteq P$ (which respectively
equals $P$ or extends it at one or both of its end-points according
as $\bs{\chi}\in C_{1,1}\cup C_{1,\infty}$, $C_{1,2}\cup
C_{2,\infty}$ or $C_{2,2}$).

Thus we have obtained the Borel partition
\[\bbZ_2^{\bfF_2} = C_0 \cup C_{1,1} \cup C_{1,2}\cup C_{2,2}\cup C_{1,\infty}\cup C_{2,\infty}\cup C_{\infty,\infty}.\]
In fact, it is easy to refine this partition even further. If
$\bs{\chi}\in C_{a,b}$ with $a,b < \infty$, then $P$ and $R$ are
finite subsets of $\bfF_2$. Moreover, the fact that $\bs{\chi}\in
C_{a,b}$ now depends only on the restriction $\bs{\chi}|_{B(R,10)}$
(in the sense that any other $\bs{\chi}'$ agreeing with $\bs{\chi}$
on this restriction also lies in $C_{a,b}$, with walkers seeing just
the same configurations).  We may therefore partition $C_{a,b}$
according to the triples $(P,R,\psi)$, where $\psi :=
\bs{\chi}|_{B(R,10)\setminus R}$, that can arise in this way.

Let $\O_{a,b}$ be the collection of triples $(P,R,\psi)$ such that
any point $\bs{\chi}$ giving rise to them as above must lie in
$C_{a,b}$, and let
\[C_{P,R,\psi} := \{\bs{\chi}\in \bbZ_2^{\bfF_2}:\ \bs{\chi}^{-1}\{0\} \supseteq R\ \hbox{and}\ \bs{\chi}|_{B(R,10)\setminus R} = \psi\}\]
be the cylinder set associated to this triple. In this situation we
will refer to $P$ as the \textbf{inner path} and $R$ as the
\textbf{outer path} of $(P,R,\psi)$.  Clearly $R = P$ if and only if
$a = b = 1$, and sometimes we will abusively write members of $\O_{1,1}$
simply as pairs $(P,\psi)$.

We have now obtained the following finer partition.

\begin{lem}\label{lem:partition}
The equality
\[\bbZ_2^{\bfF_2} = C_0\cup \Big(\bigcup_{\scriptsize{\begin{array}{c}a,b \in \{1,2\},\\ a\leq b\end{array}}}\bigcup_{(P,R,\psi)\in \O_{a,b}}C_{P,R,\psi}\Big) \cup C_{1,\infty}\cup C_{2,\infty}\cup C_{\infty,\infty}\]
holds, and is a Borel partition of $\bbZ_2^{\bfF_2}$. \qed
\end{lem}

From this partition we can obtain a related orthogonal decomposition
of the Hilbert space $L^2(m_{V^\perp}\otimes \#_{\bfF_2})$, and it is
in this form that its importance will become clear.  We will later obtain a simple description of $Q$ in terms of its behaviour
on each of these subspaces that will then enable us to identify
certain of its eigenspaces exactly. For each $(P,R,\psi) \in
\O_{a,b}$ we define
\[\frH_{P,R,\psi} := \rm{img}(M_{C_{P,R,\psi}})\]
and also
\[\frH_0 := \rm{img}(M_{C_0})\quad\quad\hbox{and}\quad\quad\frH_{a,\infty} := \rm{img}(M_{C_{a,\infty}})\quad\quad\hbox{for}\ a\in\{1,2,\infty\},\]
and so now we can write
\begin{multline*}
L^2(m_{V^\perp}\otimes \#_{\bfF_2})\\ = \frH_0 \oplus
\Big(\bigoplus_{\scriptsize{\begin{array}{c}a,b \in \{1,2\},\\a \leq
b\end{array}}}\bigoplus_{(P,R,\psi) \in
\O_{a,b}}\frH_{P,R,\psi}\Big)\oplus \frH_{1,\infty}\oplus
\frH_{2,\infty}\oplus \frH_{\infty,\infty}.
\end{multline*}

Note that since each component of this decomposition is defined by
an orthogonal projection lying in the von Neumann algebra
$L^\infty(m_{V^\perp})\rtimes \bfF_2 \cong L\G$, each defines a
submodule for the right action of $\G$ on $L^2(m_{V^\perp}\otimes
\#_{\bfF_2})$ (arising by applying the Fourier transform to the
right von Neumann algebra $R\G$ acting on $\ell^2(\G)$, as in
Section~\ref{sec:prelims}), and has a well-defined von Neumann
dimension given by the standard trace on
$L^\infty(m_{V^\perp})\rtimes \bfF_2$.

It will turn out that for a suitable choice of $V^\perp$ we have
\[m_{V^\perp}(C_{a,\infty}) = 0\quad\quad\forall a \in \{1,2,\infty\},\]
so that the spaces $\frH_{a,\infty}$ for $a \in \{1,2,\infty\}$
contribute trivially to the above decomposition.  This will be
proved in Proposition~\ref{prop:no-infinite-paths} once we have
specified our method for choosing $V$. In the remainder of this
section we make a closer examination of the behaviour of $Q$ on the
subspaces $\frH_0$ and $\frH_{P,R,\psi}$.

We first organize the above orthogonal decomposition by `clustering'
the subspaces involved into certain equivalence classes, in such a
way that the subspaces of the coarser decomposition that results
from this clustering are individually $Q$-invariant and each admits
a relatively simple description of the action of $Q$.  The
equivalence relation we need is the following.

\begin{dfn}[Translation equivalence]
Two triples $(P_1,R_1,\psi_1)$ and $(P_2,R_2,\psi_2)$, with
$P_1,P_2$ (finite or infinite) paths in $\bfF_2$ that pass within
distance $1$ of $e$ and $\psi_i:B(R_i,10)\setminus R_i \to \bbZ_2$,
are \textbf{translation equivalent} (denoted by
$(P_1,R_1,\psi_1)\sim (P_2,R_2,\psi_2)$) if there is some $g \in
\bfF_2$ such that $P_2 = gP_1$, $R_2 = gR_1$ and $\psi_2(gh) =
\psi_1(h)$ for all $h \in B(R_1,10)\setminus R_1$. In this case we
will also write that $(P_1,R_1,\psi_1)$ is a \textbf{translate} of
$(P_2,R_2,\psi_2)$. Since $P_1$ and $P_2$ are both required to pass
within distance $1$ of $e$, if $P_1$ is finite then clearly the
equivalence class of $(P_1,R_1,\psi_1)$ is a finite set of size
$|B(P_1,1)|$.
\end{dfn}

We use this to re-organize the above orthogonal decomposition as
\[L^2(m_{V^\perp}\otimes \#_{\bfF_2}) = \frH_0\oplus\Big(\bigoplus_{\scriptsize{\begin{array}{c}a,b\in\{1,2\}\\ a\leq b\end{array}}}\bigoplus_{\C \in \O_{a,b}/\sim}\frH_\C\Big)\oplus \frH_{1,\infty}\oplus \frH_{2,\infty}\oplus \frH_{\infty,\infty},\]
where
\[\frH_\C := \bigoplus_{(P,R,\psi) \in \C}\frH_{P,R,\psi}.\]

The following is a straightforward extension of Equation (3.5) in
Dicks and Schick~\cite{DicSch02}.

\begin{lem}\label{lem:Q-acting-on-comp-ops}
For any $g \in \bfF_2$ and any Borel subset $Y\subseteq
\bbZ_2^{\bfF_2}$ we have
\[Q\circ T^{g^{-1}}\circ M_Y = \sum_{s \in S}T^{s^{-1}g^{-1}}\circ (M_{(F_s\circ \hat{\a}^{g^{-1}})\cdot 1_Y} + M_{(G_s\circ \hat{\a}^{s^{-1}g^{-1}})\cdot 1_Y}).\]
\qed
\end{lem}

\begin{lem}
We have $Q|_{\frH_0} = 0$.
\end{lem}

\textbf{Proof}\quad By the definition of $C_0$ and
Lemma~\ref{lem:Q-acting-on-comp-ops} we have
\[F_s\cdot 1_{C_0} = (G_s\circ\hat{\a}^{s^{-1}})\cdot 1_{C_0} = 0\quad\quad\forall s \in S\]
and so
\begin{multline*}
M_{C_0}f = f\\ \Rightarrow\quad\quad Qf = (Q\circ M_{C_0})f
= \sum_{s\in S}\big(T^{s^{-1}}\circ (M_{F_s\cdot 1_{C_0}} +
M_{(G_s\circ\hat{\a}^{s^{-1}})\cdot 1_{C_0}})\big)f = 0.
\end{multline*}
\qed

\begin{figure}
\begin{center}
\includegraphics[width = 1.0\textwidth]{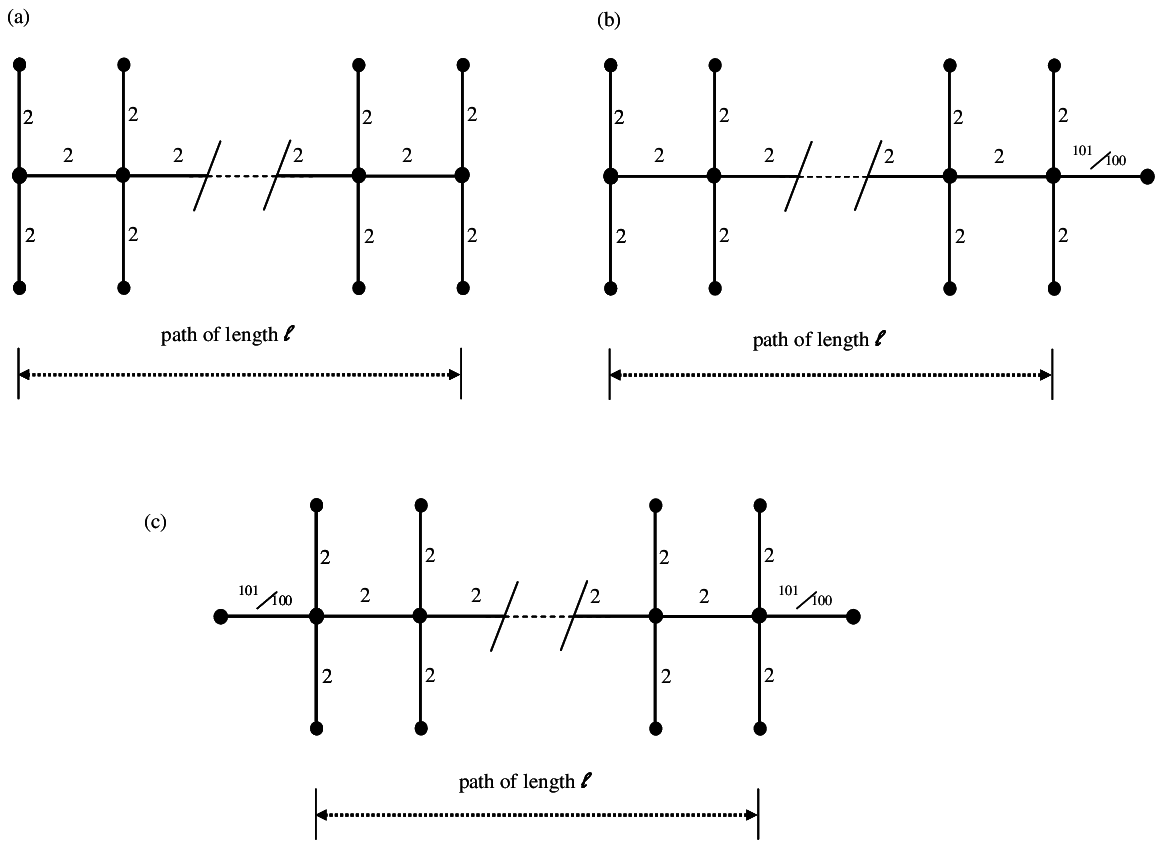}
\caption{The weighted graph $(V^{\ell,a,b},E^{\ell,a,b,})$
corresponding to $Q|_{\frH_\C}$ for (a) $\C \in \O_{1,1}/\sim$, (b)
$\C \in \O_{1,2}/\sim$ and (c) $\C \in
\O_{2,2}/\sim$}\label{fig:Qfin}
\end{center}
\end{figure}

\begin{prop}\label{prop:description-of-Q-on-fin-bits}
Let $(V^{\ell,a,b},E^{\ell,a,b})$ for $a,b \in \{1,2\}$ be the
weighted graphs shown in Figure~\ref{fig:Qfin} and
\[Q^{\ell,a,b} = (q_{u,v}^{\ell,a,b})_{u,v \in V^{\ell,a,b}}\]
their weighted adjacency matrices, regarded as operators on
$\ell^2(V^{\ell,a,b})$.  Then
for each $\C \in \O_{a,b}/\sim$ such that $(P,R,\psi)\in \C$ has
$|P| = \ell$, the subspace $\frH_\C$ is $Q$-invariant, and there is
some von Neumann right-module $\frh_\C$ (which will in fact depend
on the measure $m_{V^\perp}$) such that we have
\[Q|_{\frH_\C} \cong \id_{\frh_\C}\otimes Q^{\ell,a,b}.\]
\end{prop}

\textbf{Proof}\quad We treat the case of $Q|_{\frH_\C}$ for some
$\C\in \O_{1,2}/\sim$, the others being similar.

Pick a representative $(P,R,\psi) \in \C$, say with $|P| = \ell$,
such that $e$ is the `good' end-point of $P$: that is, such that
$\bs{\chi}$ itself is locally good. There is exactly one such
end-point if $(P,R,\psi)\in \O_{1,2}$. Let $\frh_\C$ be the von
Neumann right-module $\frH_{P,R,\psi}$ (of course, the dimension of
this depends on $m_{V^\perp}$).  Owing to the involvement of $W$ in
the definition of $F_s$ and hence of $\O_{1,2}$, we know that $\ell
\geq 5$.

Next observe that if $g \in B(P,1)$, then the values
\[F_s(\hat{\a}^{g^{-1}}(\bs{\chi}))\quad\hbox{and}\quad G_s(\hat{\a}^{s^{-1}}\hat{\a}^{g^{-1}}(\bs{\chi})) = F_{s^{-1}}(\hat{\a}^{(gs)^{-1}}(\bs{\chi}))\]
are the same for all $\bs{\chi} \in C_{P,R,\psi}$.  In view of this
we can define
\[\phi(g,gs) := F_s(\hat{\a}^{g^{-1}}(\bs{\chi}))\]
for $g \in B(P,1)$ using any representative $\bs{\chi}\in
C_{P,R,\psi}$, and obtain
\begin{multline}\label{eq:helpful}
M_{(F_s\circ \hat{\a}^{g^{-1}})\cdot 1_{C_{P,R,\psi}}} +
M_{(G_s\circ \hat{\a}^{s^{-1}g^{-1}})\cdot 1_{C_{P,R,\psi}}}\\ =
(\phi(g,gs) + \phi(gs,g))\cdot M_{C_{P,R,\psi}}.
\end{multline}

We can now simply read off from the definition of $F_s$ a very
explicit description of this function $\phi$ on the set of pairs
\[\{(g,h):\ g,h \in B(P,1),\,\rho(g,h) = 1\}:\]
\begin{itemize}
\item If $g$ is an interior point of $P$, then it has
\begin{itemize}
\item two neighbours $h$ that are not in $P$, and for each of these we
have $\phi(g,h) = 2$ and $\phi(h,g) = 0$, so $\phi(g,h) + \phi(h,g)
= 2$, and
\item two neighbours $h$ that are also in $P$, so if such an $h$ is
also an interior point then $\phi(g,h) = \phi(h,g) = 1$ and if it is
an end-point of $P$ then $\phi(g,h) = 2$ and $\phi(h,g) = 0$, and in
either case overall $\phi(g,h) + \phi(h,g) = 2$;
\end{itemize}
\item If $g = e$ is the good end-point, then it has
\begin{itemize}
\item one neighbour $s$ that must lie in the interior of $P$, for which $\phi(e,s) = 0$ and $\phi(s,e) =
2$ and so $\phi(e,s) + \phi(s,e) = 2$,
\item an opposite neighbour $s^{-1}$, for which $\phi(e,s^{-1}) = \phi(s^{-1},e) =
0$, so $\phi(e,s^{-1}) + \phi(s^{-1},e) = 0$, and
\item two neighbours $t$ neither of which lie in $P$ and such that
$e$ is their mid-point, for each of which $\phi(e,t) = 2$ and
$\phi(t,e) = 0$ so that $\phi(e,t) + \phi(t,e) = 2$;
\end{itemize}
\item If $g$ is the other (`bad') end-point of $P$, so that it still lies in the interior of $R$, then it has
\begin{itemize}
\item one neighbour $h$ that lies in the interior of $P$, for which $\phi(g,h) = 1$ and $\phi(h,g) =
1$ and so $\phi(g,h) + \phi(h,g) = 2$,
\item one neighbour $h$ that lies in $R\setminus P$, for which $\phi(g,h) = 1$ and $\phi(h,g) =
\frac{1}{100}$, so $\phi(g,h) + \phi(h,g) = \frac{101}{100}$, and
\item two neighbours $h$ which do not lie in $R$, for each of which $\phi(g,h) = 2$ and
$\phi(h,g) = 0$ so that $\phi(g,h) + \phi(h,g) = 2$.
\end{itemize}
\end{itemize}
Note that the cases above involving the `good' end-point are where
we have used the quirk in the definition of $F_s$ discussed in
Remark 2 after Definition~\ref{dfn:FScirc}.

Putting these possibilities together, and comparing them with
Figure~\ref{fig:Qfin}(b), we see that if we let $V_0 \subseteq
V^{\ell,1,2}$ be the subset of $\ell$ vertices on the central path
of that graph then we may choose a bijection $\xi_0:V_0 \to P$ such
that the left (respectively right) end-point of $V_0$ is sent to $e$
(respectively, to the `bad' end-point of $P$), and now extend this
to an isomorphism of weighted graphs
\begin{multline*}
\xi:(V^{\ell,1,2},E^{\ell,1,2},Q^{\ell,1,2})\\ \to
\big(B(P,1),\Cay(\bfF_2,S)|_{B(P,1)},(\phi(g,h) + \phi(h,g))_{g,h
\in B(P,1),\rho(g,h) = 1}\big),
\end{multline*}
(where we have been just a little sloppy, in that we allow our
`isomorphism of weighted graphs' to miss the isolated neighbour of
$e$ with no positive-weight connections).  That this is possible
follows by inspection of Figure~\ref{fig:Qfin}(b) and the list of
possibilities above, which shows that for each $v \in V_0$ we may
pair up its neighbours with those of $\xi_0(v) \in P$ so as to
respect the edge-weights:
\[\{u,v\} \in E^{\ell,1,2}\quad\quad\Rightarrow\quad\quad\phi(\xi(u),\xi(v)) + \phi(\xi(v),\xi(u)) = q^{\ell,1,2}_{u,v}.\]

We can now simply turn this isomorphism of weighted graphs into an
isomorphism of Hilbert space operators as follows. Let
$(\delta_v)_{v \in V^{\ell,1,2}}$ be the standard basis of
$\ell^2(V^{\ell,1,2})$.  Observe from the definition of translation
equivalence that
\[\C = \{(gP,gR,\psi(g\,\cdot\,)):\ g\in B(P,1)\}\]
and that $C_{gP,gR,\psi(g\,\cdot\,)} =
\hat{\a}^{g^{-1}}(C_{P,R,\psi})$, and hence that
\begin{multline*}
\frH_\C = \bigoplus_{(P',R',\psi')\sim (P,R,\psi)}\frH_{P',R',\psi'}
= \bigoplus_{g \in
B(P,1)}\rm{img}\big(M_{\hat{\a}^{g^{-1}}(C_{P,R,\psi})}\big)\\ =
\bigoplus_{g \in B(P,1)}\rm{img}\big(T^{g^{-1}}\circ
M_{C_{P,R,\psi}}\circ T^g\big) = \bigoplus_{g \in
B(P,1)}T^{g^{-1}}(\frh_\C)
\end{multline*}

Now define
\[\Phi:\frH_\C \to \frh_\C \otimes \ell^2(V^{\ell,1,2})\]
by setting
\[\Phi(f) = T^g(f)\otimes \delta_{\xi^{-1}(g)}\quad\quad\hbox{for}\ g \in B(P,1),\ f \in T^{g^{-1}}(\frh_\C).\]

This is clearly an isomorphism of von Neumann right-modules, and it
is now simple to check that $Q|_{\frH_\C} = \Phi^{-1}\circ
(\id_{\frh_\C}\otimes Q^{\ell,1,2})\circ \Phi$: indeed, if $f \in
T^{g^{-1}}(\frh_\C)$, so $M_{C_{P,R,\psi}}(T^g f) = T^g f$, then
using Lemma~\ref{lem:Q-acting-on-comp-ops} and
equation~(\ref{eq:helpful}) we have
\begin{eqnarray*}
Qf &=& (Q\circ T^{g^{-1}})(T^g f) = (Q\circ
T^{g^{-1}})(M_{C_{P,R,\psi}}(T^g f))\\
&=& (Q\circ T^{g^{-1}}\circ M_{C_{P,R,\psi}})(T^gf)\\
&=& \sum_{s\in S}(\phi(g,gs) +
\phi(gs,g))\cdot T^{s^{-1}g^{-1}}\circ M_{C_{P,R,\psi}}(T^gf)\\
&=& \sum_{s\in S}(\phi(g,gs) +
\phi(gs,g))\cdot T^{s^{-1}g^{-1}}(T^gf)\\
&=& \sum_{s\in S}T^{s^{-1}g^{-1}}(q^{\ell,1,2}_{\xi^{-1}(g),\xi^{-1}(gs)}(T^gf))\\
&=& \Phi^{-1}\Big(T^gf\otimes \Big(\sum_{s\in S}q^{\ell,1,2}_{\xi^{-1}(g),\xi^{-1}(gs)}\delta_{\xi^{-1}(gs)}\Big)\Big)\\
&=& \Phi^{-1}(T^gf\otimes (Q^{\ell,1,2}(\delta_{\xi^{-1}(g)}))) =
(\Phi^{-1}\circ(\id_{\frh_\C}\otimes Q^{\ell,1,2})\circ\Phi)(f),
\end{eqnarray*}
as required. \qed

\section{Computation of an eigenspace}

We now specialize to a particular selection of an eigenvalue of
interest to us: the value $4$.  We will find (in
Corollary~\ref{cor:eigenspace-for-4} below) that we can describe the
eigenspace $\ker(Q - 4\cdot \id)$ rather explicitly.  The choice of
$4$ here is important: it originates in a particular quadratic
equation that arises from a two-step linear recursion that will
appear repeatedly below in the description of the associated
eigenspaces.  For this value we can obtain a proof of the
non-existence of such eigenspaces for some of the restrictions
$Q|_{\frH_\C}$, and an explicit construction of these eigenspaces
for others.

We will find that (in much the same way as for the simple
lamplighter group as described in Dicks and Schick~\cite{DicSch02})
we can arrange a `pileup' of infinitely many eigenspaces
corresponding to this eigenvalue, each of them admitting a
relatively simple description, and it is this that will ultimately
give us the control over von Neumann dimensions required for
Theorem~\ref{thm:main}.

\begin{lem}\label{lem:sometimes-excluding-4}
The value $4$ is not an eigenvalue of $Q|_{\frH_\C}$ for any $\C\in
\O_{1,2}/\sim$ or $\C \in \O_{2,2}/\sim$.
\end{lem}

\textbf{Proof}\quad We give the proof for $\C \in \O_{1,2}/\sim$,
the other case being exactly similar.  By
Proposition~\ref{prop:description-of-Q-on-fin-bits} it will suffice
to show that
\[\hbox{for}\ \ell \geq 5,\ \bf{x} \in \ell^2(V^{\ell,1,2}),\quad\quad Q^{\ell,1,2}\bf{x} = 4\bf{x}\quad\Rightarrow\quad \bf{x} = \bs{0}.\]

To this end, enumerate the central length-$(\ell-1)$ path of
$V^{\ell,1,2}$ as $V_0 = \{v_1,v_2,\ldots,v_\ell\}$, and observe
that each pair of neighbours in this path is joined by an edge of
weight $2$, that each vertex in this path is joined to exactly two
non-members of this path by edges of weight $2$, and that one of the
end-points is additionally joined to a non-member of this path by an
edge of weight $\frac{101}{100}$. These are all the positive-weight
edges in the graph. For each $i \in \{1,2,\ldots,\ell-1\}$ let
$v_{i,j}$ for $j=1,2$ be the two neighbours of $v_i$ in
$V^{\ell,1,2}\setminus V_0$, and also let $v_{\ell,1}$ and
$v_{\ell,2}$ be the two outer neighbours of $v_\ell$ joined to it by
a weight of $2$, and $v_\ell'$ the neighbour of $v_\ell$ joined to
it by the weight $\frac{101}{100}$. Finally let
\[\o := \frac{1 + \sqrt{-3}}{2},\]
so that $\{1,\o,\o^2,-1,\bar{\o}^2,\bar{\o}\}$ are the sixth roots
of unity.

We can evaluate the equation $Q^{\ell,1,2}\bf{x} = 4\bf{x}$ at $v_i$
for $2 \leq i \leq \ell - 1$ and also at $v_{i,j}$ for such $i$, and
find that
\begin{eqnarray*}
\hbox{at}\ v_{i,j}&:&\quad\quad 4\bf{x}(v_{i,j}) =
2\bf{x}(v_i)\quad\Rightarrow\quad \bf{x}(v_{i,j}) =
\frac{1}{2}\bf{x}(v_i)\\
\hbox{at}\ v_i&:&\quad\quad 4\bf{x}(v_i) = 2\bf{x}(v_{i-1}) +
2\bf{x}(v_{i+1}) + 2\bf{x}(v_{i,1}) + 2\bf{x}(v_{i,2})\\
&&\quad\quad\quad\Rightarrow\quad \bf{x}(v_i) = \bf{x}(v_{i-1}) +
\bf{x}(v_{i+1}),
\end{eqnarray*}
and so re-arranging we obtain
\[\bf{x}(v_{i+1}) = \bf{x}(v_i) - \bf{x}(v_{i-1})\quad\quad \forall i=2,3,\ldots,\ell - 1,\]
and hence by solving this quadratic recursion that there are $a,b
\in \bbC$ such that $\bf{x}(v_i) = a\o^i + b\bar{\o}^i$; and we also
obtain similarly that $\bf{x}(v_{i,j}) = \frac{1}{2}\bf{x}(v_i)$ for
all $i\in \{1,2,\ldots,\ell\}$ and $j = 1,2$.

Next, evaluating at $v_\ell'$ gives that $\bf{x}(v_\ell') =
\frac{101}{400}\bf{x}(v_\ell)$, and now evaluating at $v_\ell$ gives
\begin{eqnarray*}
4\bf{x}(v_\ell) &=& 2\bf{x}(v_{\ell-1}) + 2(\bf{x}(v_{\ell,1}) +
\bf{x}(v_{\ell,2})) +
\frac{101^2}{4\cdot 100^2}\bf{x}(v_\ell)\\
&=& 2\bf{x}(v_{\ell - 1}) + 2\bf{x}(v_\ell) + \frac{101^2}{4\cdot
100^2}\bf{x}(v_\ell)\\
&\Rightarrow&\quad\quad \Big(1 - \frac{101^2}{8\cdot
100^2}\Big)(a\o^\ell + b\bar{\o}^\ell) = a\o^{\ell - 1} +
b\bar{\o}^{\ell - 1}.
\end{eqnarray*}
It follows that either $\bf{x} = \bs{0}$ or at least one of $a,b$ is
non-zero.  Let us suppose it is $b$ and derive a contradiction, the
case $a \neq 0$ being similar.  In this case the above conclusion
can be re-arranged to give
\[\Big(\Big(1 - \frac{101^2}{8\cdot 100^2}\Big)\o^{\ell} - \o^{\ell - 1}\Big)\frac{a}{b} = \bar{\o}^{\ell-1} - \Big(1 - \frac{101^2}{8\cdot 100^2}\Big)\bar{\o}^\ell,\]
and now evaluating the eigenvector equation at $v_1$ (the only
vertex where we have not yet checked it) gives similarly \[a\o +
b\bar{\o} = \bf{x}(v_1) = \bf{x}(v_2) = a\o^2 + b\bar{\o}^2
\quad\quad\Rightarrow\quad\quad (\o - \o^2)\frac{a}{b} = \bar{\o}^2
- \bar{\o}.\]

It can now be verified directly that no value $\frac{a}{b}$ can
simultaneously satisfy both of the above equations (bearing in mind
that the sequence $(\o^\ell)_{\ell \geq 1}$ takes only six values).
This gives the desired contradiction, and so completes the proof.
\qed

\begin{lem}
The value $4$ is an eigenvalue of $Q^{\ell,1,1}$ with multiplicity
$1$ whenever $\ell \equiv -1 \mod 6$, and hence also of
$Q|_{\frH_\C}$ for any $\C \in \O_{1,1}/\sim$ such that $(P,\psi)\in
\C$ has $|P|\equiv -1 \mod 6$ and $\frh_\C \neq \{0\}$.
\end{lem}

\textbf{Proof}\quad If $\frh_\C \neq \{0\}$ then the conclusion for
$Q|_{\frH_\C}$ follows directly from that for $Q^{\ell,1,1}$ using
Proposition~\ref{prop:description-of-Q-on-fin-bits}, so we focus on
the latter.  We explicitly exhibit a suitable eigenvector, and then
the argument of the preceding lemma shows that it is the only one up
to scalar multiples. Let $V_0 = \{v_1,v_2,\ldots,v_\ell\}$,
$v_{i,j}$ for $j=1,2$ and $\o \in\bbC$ be as in the preceding lemma,
and now define $\bf{x} \in \ell^2(V^{\ell,1,1})$ by
\[\bf{x}(v_i) := \o^i + \bar{\o}^{3+i}\quad\quad\hbox{and}\quad\quad \bf{x}(v_{i,j}) := \frac{1}{2}(\o^i + \bar{\o}^{3+i})\]
for $v_i$, $v_{i,j} \in V^{\ell,1,1}$.

It is now a simple check that $Q^{\ell,1,1}\bf{x} = 4\bf{x}$:
\begin{itemize}
\item At $v_i$ for $2 \leq i\leq \ell - 1$ we have
\begin{eqnarray*}
(Q^{\ell,1,1}\bf{x})(v_i) &=& 2\bf{x}(v_{i-1}) + 2\bf{x}(v_{i+1}) +
2\bf{x}(v_{i,1}) + 2\bf{x}(v_{i,2})\\
&=& 2(\o^{i-1} + \bar{\o}^{3 + i-1} + \o^{i+1} + \bar{\o}^{3 + i+1}
+ \o^i +
\bar{\o}^{3 + i})\\
&=& 2\big((\o^i + \bar{\o}^{3+i}) + (\o^i + \bar{\o}^{3+i})\big)\\
&=& 4(\o^i + \bar{\o}^{3+i}) = 4\bf{x}(v_i),
\end{eqnarray*}
since $\o + \bar{\o} = 1$.
\item At $v_1$ we have $\bf{x}(v_1) = \o + \bar{\o}^4 = \o + \o^2$,
and
\begin{multline*}
(Q^{\ell,1,1}\bf{x})(v_1) = 2\bf{x}(v_2) + 2\bf{x}(v_{1,1}) +
2\bf{x}(v_{2,2})\\
= 2\Big(\o^2 + \bar{\o}^5 + 2\cdot \frac{1}{2}(\o + \o^2)\Big) =
4(\o + \o^2) = 4\bf{x}(v_1),
\end{multline*}
using now that $\bar{\o}^5 = \o$, and similarly since $\ell \equiv
-1 \mod 6$ we have
\begin{multline*}
(Q^{\ell,1,1}\bf{x})(v_\ell) = 2\bf{x}(v_{\ell - 1}) +
2\bf{x}(v_{\ell,1}) +
2\bf{x}(v_{\ell,2})\\
= 2\Big(\o^{\ell - 1} + \bar{\o}^{\ell + 2} + 2\cdot
\frac{1}{2}(\o^\ell + \bar{\o}^{3 + \ell})\Big) = 4(\o^{-1} +
\bar{\o}^2) = 4\bf{x}(v_\ell)
\end{multline*}
(notice, however, that more general linear combinations of $\o^i$
and $\bar{\o}^i$ would not work here);
\item Finally, at a leaf $v_{i,j}$ we need $2\bf{x}(v_i) =
4\bf{x}(v_{i,j})$, and this is obvious. \qed
\end{itemize}

Combining the above calculations now gives the following.

\begin{cor}\label{cor:eigenspace-for-4}
With $\G := (\bbZ_2^{\oplus \bfF_2}/V)\rtimes \bfF_2$, under the
assumption that $\frH_{1,\infty} = \frH_{2,\infty} =
\frH_{\infty,\infty} = \{0\}$ we have
\[\ker(Q - 4\cdot \id_{L^2(m_{V^\perp}\otimes \#_{\bfF_2})}) = \bigoplus_{i\geq 1}\ker(Q|_{\frH_{\C_i}} - 4\cdot\id_{\frH_{\C_i}})\]
for some infinite sequence $\C_1$, $\C_2$, \ldots, in
$\O_{1,1}/\sim$, and hence
\begin{eqnarray*}
\dim_{L\G}\ker(Q - 4\cdot \id_{L_2(m_{V^\perp}\otimes \#_{\bfF_2})})
&=& \sum_{i\geq 1}\dim_{L\G}\ker(Q|_{\frH_{\C_i}} - 4\cdot
\id_{\frH_{\C_i}})\\
&=& \sum_{i\geq 1}\dim_{L\G}\,\frh_{\C_i}.
\end{eqnarray*}
\end{cor}

\textbf{Proof}\quad This all follows directly from the preceding
lemmas upon noting that since the value $4$ has multiplicity $1$ as
an eigenvalue of $Q^{\ell,1,1}$ we have
\begin{multline*}
\dim_{L\G}\ker(Q|_{\frH_{\C_i}} - 4\cdot \id_{\frH_{\C_i}}) =
\dim_{L\G}\ker(\id_{\frh_{\C_i}}\otimes (Q^{\ell,1,1}\ - 4\cdot
\id_{\ell^2(V^{\ell,1,1})})\\ = \dim_{L\G}\,\frh_{\C_i}.
\end{multline*}
\qed

\begin{dfn}
We will refer to those $\C \in \O_{1,1}/\sim$ that contribute
nontrivially to the above sum expression for $\dim_{L\G}\ker(Q -
4\cdot \id_{L_2(m_{V^\perp}\otimes \#_{\bfF_2})})$, or also any
$(P,\psi)$ that lies in such a $\C$, as \textbf{active}.
\end{dfn}

\section{Estimates on von Neumann dimensions}

So far our results have been independent of the particular choice of
the subspace $V$, and in particular of the Haar measure
$m_{V^\perp}$, even though it has already been mentioned in the
notation a number of times. That choice will now become important,
as we seek to show how certain possible choices of $V$ give
different possible values for the von Neumann dimensions of the
subspaces in Corollary~\ref{cor:eigenspace-for-4}.

The calculation of these dimensions will rest on the following
lemma.

\begin{lem}\label{lem:measure-equipartition}
Suppose that $V \leq \bbZ_2^{\bfF_2}$ is a subgroup, $A \subset
\bfF_2$ is a finite subset and for $\psi:A\to\bbZ_2$ let
\[C(\phi) := \{\bs{\chi}\in\bbZ_2^{\bfF_2}:\ \bs{\chi}|_A = \phi\}.\]
Then
\[m_{V^\perp}(C(\phi)) = \left\{\begin{array}{ll}\frac{1}{|\{\phi' \in \bbZ_2^A:\ C(\phi')\cap V^\perp \neq \emptyset\}|}&\quad \hbox{if}\ C(\phi)\cap V^\perp \neq \emptyset\\ 0&\quad\hbox{else}\end{array}\right.\]
(that is, the measure $m_{V^\perp}$ is shared equally among those
cylinder sets $C(\phi)$ that intersect $V^\perp$ nontrivially).
\end{lem}

\textbf{Proof}\quad Clearly $m_{V^\perp}(C(\phi)) = 0$ if
$C(\phi)\cap V^\perp = \emptyset$, so it suffices to prove that
every $C(\phi)$ for which $C(\phi)\cap V^\perp \neq \emptyset$ has
equal measure under $m_{V^\perp}$. If $C(\phi_1),C(\phi_2)$ are two
such, then we can pick some $\bs{\chi}_i \in C(\phi_i)\cap V^\perp$
for $i = 1,2$, and now inside the group $\bbZ_2^{\bfF_2}$ the
translation by $\bs{\chi}_2 - \bs{\chi}_1$ is
$m_{V^\perp}$-preserving and sends $C(\phi_1)\cap V^\perp$ to
$C(\phi_2)\cap V^\perp$, so this completes the proof. \qed

We now turn to the steps needed in our construction of the subgroups
$V_I$.  Our first step is to pick a strictly increasing sequence
$(l(n))_{n\geq 1}$ in $\bbN$ (where we adopt the convention
$0\not\in\bbN$).

\begin{lem}\label{lem:free1}
The elements $t_i := s_2^{l(i)} s_1 s_2^{-l(i)}$, $i\geq 1$, are
free in $\bfF_2$, and so generate a homomorphic embedding
$\bfF_\infty \into \bfF_2$.
\end{lem}

\textbf{Proof}\quad Suppose that
\[t_{i_1}^{k_1}t_{i_2}^{k_2}\cdots t_{i_m}^{k_m} = e\]
for some sequences $i_1$, $i_2$, \ldots, $i_m \in \{0,1,\ldots,n\}$
and $k_1$, $k_2$, \ldots, $k_m \in\bbZ\setminus \{0\}$. Then since
$t_i^k = s_2^{l(i)} s_1^{k} s_2^{-l(i)}$ for all $k\in \bbZ$, we may
reduce this evaluation to
\[s_2^{l(i_1)} s_1^{k_1} s_2^{l(i_2) - l(i_1)} s_1^{k_2} s_2^{l(i_3) - l(i_2)} \cdots s_2^{k_m} s_1^{-l(i_m)} = e,\]
and it is now clear that this is possible only if $i_1 = i_2 =
\ldots = i_m$ and $k_1 + k_2 + \cdots + k_m = 0$, hence only if the
original word was trivial. \qed

\begin{lem}\label{lem:free2}
For any $h \in \langle t_n:\ n\geq 1\rangle\setminus \{e\}$ the path
in $\Cay(\bfF_2,S)$ joining $e$ to $h$ passes through $s_2^{\pm 1}$
and not through $s_1^{\pm 1}$.
\end{lem}

\textbf{Proof}\quad If \[h = t_{i_1}^{k_1}t_{i_2}^{k_2}\cdots
t_{i_m}^{k_m}\] for some $i_1$, $i_2$, \ldots, $i_m \in
\{1,2,\ldots,n\}$ with consecutive values distinct and some $k_1$,
$k_2$, \ldots, $k_m \in \bbZ \setminus \{0\}$, then as before we can
write this out as
\[s_2^{l(i_1)} s_1^{k_1} s_2^{l(i_2) - l(i_1)} s_1^{k_2} s_2^{l(i_3) - l(i_2)} \cdots s_2^{k_m} s_1^{-l(i_m)},\]
and this is now the reduced word form of $h$. Since the path in
question is just the sequence of initial segments of this word, we
can see that the first step must be $s_2^{\pm 1}$, as required. \qed

Now for $I \subseteq \bbN$ we define
\[V_I := \rm{span}_{\bbZ_2}\Big\{\sum_{i=-10}^{10}(\delta_{gs_1^i} - \delta_{gt_ns_1^i}):\ g \in \bfF_2,\ n\in I\Big\},\]
so that
\[V_I^\perp := \Big\{\bs{\chi}\in \bbZ_2^{\bfF_2}:\ \sum_{i=-10}^{10}\bs{\chi}(gs_1^i) = \sum_{i=-10}^{10}\bs{\chi}(gt_ns_1^i)\ \forall g\in \bfF_2,\,n\in I\Big\}.\]

Let us also write $\G_I := (\bbZ_2^{\oplus \bfF_2}/V_I)\rtimes
\bfF_2$ and let $Q_I$ be the operator in $\bbQ\G_I$ defined as in
Section~\ref{sec:intro-operators}.  Finally, let
\[\L_I := \langle t_n:\ n\in I\rangle \leq \bfF_2,\]
and for any subset $A \subseteq \bfF_2$ let
\[A/\L_I := \{A\cap g\L_I:\ g\in A\},\]
the partition of $A$ induced by the partition of $\bfF_2$ into
left-cosets of $\L_I$.

\begin{lem}\label{lem:dim-formula}
Let $\C_1$, $\C_2$, \ldots $\in \O_{1,1}/\sim$ be the active
equivalence classes, and for each $i \in \bbN$ let
$(P_i,\psi_i)\in\C_i$ be a representative for which $\frh_{\C_i} =
\frH_{P_i,\psi_i}$. Then for any $I \subseteq \bbN$ we have
\begin{eqnarray}
\dim_{L\G_I}\ker(Q_I - 4) &=& \sum_{i\geq
1}m_{V_I^\perp}(C_{P_i,\psi_i})\\ &=& \sum_{i\geq
1}\frac{1_{\{C_{P_i,\psi_i}\cap V_I^\perp\neq \emptyset\}}}{|\{\phi
\in \bbZ_2^{B(P_i,10)}:\ C(\phi)\cap V_I^\perp
 \neq \emptyset\}|}.
\end{eqnarray}
\end{lem}

\textbf{Proof}\quad This follows simply from evaluating the
individual terms in the right-hand side of
Corollary~\ref{cor:eigenspace-for-4} and observing directly from the
formula for the trace on $L^\infty(m_{\hat{U}})\rtimes \bfF_2$ that
\[\dim_{L\G_I}\,\frh_{\C_i} = \rm{tr}_{L\G_I}\,M_{C_{P_i,\psi_i}} =
m_{V^\perp}(C_{P_i,\psi_i}).\] The second line now follows from
Lemma~\ref{lem:measure-equipartition}. \qed

Next we need a criterion for deciding whether $C(\phi)\cap V_I^\perp
= \emptyset$.

\begin{lem}[Extensibility lemma]
If $A \subseteq \bfF_2$ is connected in the graph $\Cay(\bfF_2,S)$
and $\phi:A \to \bbZ_2$ is such that
\begin{eqnarray}\label{eq:local-constraints}
\sum_{i=-10}^{10}\phi(gs_1^i) = \sum_{i=-10}^{10}\phi(ghs_1^i)
\end{eqnarray}
whenever $g\in \bfF_2$ and $h\in \L_I$ are such that $gs_1^{-10}$,
$gs_1^{-9}$, \ldots, $gs_1^{10}$, $ghs_1^{-10}$, $ghs_1^{-9}$,
\ldots and $ghs_1^{10}$ all lie in $A$, then $\phi$ admits an
extension $\bs{\chi} \in C(\phi)\cap V_I^\perp$.
\end{lem}

\textbf{Remark}\quad Both the connectedness assumption on $A$ and
the fact that the removal of any vertex from $\Cay(\bfF_2,S)$
disconnects this graph are important for this proof. \fin

\textbf{Proof}\quad Given $\bs{\chi} \in C(\phi)$, it is a member of
$V_I^\perp$ if and only if for some (and hence any) upwards directed
family of subsets $B \subseteq \bfF_2$ that covers all of $\bfF_2$
we have that for each $B$ in the family the
condition~(\ref{eq:local-constraints}) holds whenever $gs_1^{-10}$,
$gs_1^{-9}$, \ldots, $gs_1^{10}$, $ghs_1^{-10}$, $ghs_1^{-9}$,
\ldots and $ghs_1^{10}$ all lie in $B$. Now let \[A_0 = A \subset
A_1\subset A_2 \subset \ldots \subset \bfF_2\] be an exhaustion of
$\bfF_2$ in which each $A_{n+1}$ is obtained from $A_n$ by the
inclusion of a single new point from $B(A_n,1)\setminus A_n$
(clearly such an exhaustion exists).  If we show how to construct
recursively a sequence of functions $\bs{\chi}_n:A_n\to\bbZ_2$ for
$n\geq 0$ such that
\begin{itemize}
\item $\bs{\chi}_0 := \phi$,
\item $\bs{\chi}_{n+1}|_{A_n} = \bs{\chi}_n$ for all $n\geq 0$
and
\item condition~(\ref{eq:local-constraints}) is satisfied by $\bs{\chi}_n$
whenever $gs_1^{-10}$, \ldots, $gs_1^{10}$, $ghs_1^{-10}$, \ldots
and $ghs_1^{10}$ all lie in $A_n$,
\end{itemize}
then it follows that $(\cup_{n\geq 1}\bs{\chi}_n) \in C(\phi)\cap
V^\perp$ is the desired point.

Moreover, having set $\bs{\chi}_0 := \phi$, it suffices to give the
construction for $\bs{\chi}_1$, since then simply repeating this
construction with $A_n$ in place of $A$ at every step completes the
proof.

To this end, suppose $A_1 = A \cup \{g_1\}$, let us write
$\cal{E}(A)$ for the set of all equations of the
form~(\ref{eq:local-constraints}) for which $gs_1^{-10}$, \ldots,
$gs_1^{10}$, $ghs_1^{-10}$, \ldots and $ghs_1^{10}$ all lie in
$A_1$, and let us partition this as
\[\cal{E}(A) = \cal{E}_0(A)\cup \cal{E}_1(A),\]
where $\cal{E}_0(A)$ contains those equations that do not involve
the value of $\bs{\chi}_1(g_1)$ and $\cal{E}_1(A)$ contains those that
do. All members of $\cal{E}_0(A)$ are satisfied by our assumptions
on $\phi$, whereas each member of $\cal{E}_1(A)$ prescribes a value
for $\bs{\chi}_1(g_1)$ in terms of values of $\phi$.  If $\cal{E}_1(A)
= \emptyset$ then we may adopt either possible value for
$\bs{\chi}_1(g_1)$, so it suffices to show that if $\cal{E}_1(A) \neq
\emptyset$ then all the resulting prescriptions agree. To see this,
observe that any two of these equations from $\cal{E}_1(A)$ must
take the form
\[\bs{\chi}(g_1) = -\sum_{-10 \leq i \leq 10,\,gs_1^i \neq g_1}\phi(gs_1^i) + \sum_{i=-10}^{10}\phi(gh_is_1^i)
\]
for some $g \in A$ and $h_1,h_2 \in \L_I$.  However, if
\[g_1 \in \{gs_1^{-10},gs_1^{-9},\ldots,gs_1^{10}\} \subset A_1,\]
then $g_1$ must be one of the end-points $gs_1^{\pm 10}$, for
otherwise $g_1 \not\in A$ would separate $A$ into the two connected
components containing these two end-points, contrary to our
assumption that $A$ is connected.  Moreover, if $g' \in A$ is
another point such that
\[g \in \{g's_1^{-10},g's_1^{-9},\ldots,g's_1^{10}\} \subset A_1,\]
then we must have $g = g'$, for if alternatively $gs_1^{-10} = g_1 =
g's_1^{10}$ then $g_1$ disconnects the components of $A$ that
contain $g$ and $g'$.  Hence we may assume without loss of
generality that all of the above equations from the collection
$\cal{E}_1(A)$ have $gs_1^{10} = g_1$.  However, since $h_1^{-1}h_2
\in \L_I$, we now see that the equation
\[\sum_{i=-10}^{10}\phi((gh_1)s_1^i) = \sum_{i=-10}^{10}\phi((gh_1)h_1^{-1}h_2s_1^i) = \sum_{i=-10}^{10}\phi(gh_2s_1^i)\]
is a member of $\cal{E}_0(A)$ and so is satisfied by assumption;
this implies that the right-hand-sides above are equal for the
equations in $\cal{E}_1(A)$ corresponding to $h_1$ and to $h_2$, and
hence prescribe a consistent value for $\bs{\chi}_1(gs_1^{10})$, as
required. \qed

\begin{cor}\label{cor:criterion-for-nonempty-cylinder}
If $P \subset \bfF_2$ is a path and $\phi:B(P,10)\to\bbZ_2$ then
$C(\phi)\cap V_I^\perp \neq \emptyset$ if and only the function
\[P\to\bbZ_2:g\mapsto \sum_{i=-10}^{10}\phi(gs_1^i)\]
is constant on the cells of $P/\L_I$.
\end{cor}

\textbf{Proof}\quad The necessity is obvious, and the sufficiency
follows from the previous lemma and the fact that
\[\{gs_1^{-10},gs_1^{-9},\ldots,gs_1^{10}\} \subset B(P,10)\quad\quad\Rightarrow\quad\quad g\in P.\]
This follows from the connectedness of $P$, because there must be
some $g_1,g_2 \in P$ that lie within distance $10$ of $gs_1^{-10}$
and $gs_1^{10}$ respectively, and were $g$ not itself a member of
$P$ then these two other members of $P$ would occupy distinct
connected components, giving a contradiction. \qed

\begin{cor}\label{cor:always-get-zero}
If $P \subset \bfF_2$ is a path with no small horizontal doglegs and
$\phi:B(P,10)\to \bbZ_2$ takes the value $0$ inside $P$ and $1$ on
$B(P,10)\setminus P$ then $C(\phi)\cap V_I^\perp \neq \emptyset$.
\end{cor}

\textbf{Remark}\quad It is in this proof that we will finally see
the purpose of the assumption of no small horizontal doglegs. \fin

\textbf{Proof}\quad By the previous corollary this depends only on
the constancy of the values
\[\sum_{i=-10}^{10}\phi(gs_1^i),\ g \in P\]
on each cell of $P/\L_I$. For $\phi$ as described this value is just
\begin{multline*}
|\{gs_1^{-10},gs_1^{-9},\ldots,gs_1^{10}\}\cap (B(P,10)\setminus
P)|\mod 2\\
\equiv|\{gs_1^{-10},gs_1^{-9},\ldots,gs_1^{10}\}\setminus P|\mod 2.
\end{multline*}
If $g$ lies in a singleton cell of $P/\L_I$ then there is nothing to
check. On the other hand, if $g,gh \in P$ for some $h \in
\L_I\setminus \{e\}$, then by applying Lemma~\ref{lem:free2} to the
segment of $P$ joining $g$ and $gh$ it follows that we must have
$gs_2^\eta \in P$ for some $\eta = \pm 1$ and $ghs_2^\eta$ for some
$\eta = \pm 1$. From this it follows that $gs_1^{\pm 1}$ cannot both
lie in $P$ and that $ghs_1^{\pm 1}$ cannot both lie in $P$.  Hence
the intersection
\[P\cap \{gs_1^{-10},gs_1^{-9},\ldots,gs_1^{10}\}\]
is either just $\{g\}$, in which case
\[|\{gs_1^{-10},gs_1^{-9},\ldots,gs_1^{10}\}\setminus P| = 20 \equiv 0\mod 2;\]
or else it also contains some point $gs_1^a$ with $a \neq 0$, so
that by the assumption of no small horizontal doglegs it must in
fact contain exactly one of the whole branches
\[\{gs_1^{-10},gs_1^{-9},\ldots,g\}\quad\quad\hbox{or}\quad\quad\{g,gs_1,\ldots,gs_1^{10}\},\]
in which case
\[|\{gs_1^{-10},gs_1^{-9},\ldots,gs_1^{10}\}\setminus P| = 10 \equiv 0\mod 2.\]

Thus the value in question is always $0 \in \bbZ_2$ for those $g$
lying in a nonsingleton cell of $P/\L_I$, and so we have proved the
necessary constancy on these cells. \qed

We will now use the preceding lemmas and corollaries to two distinct
ends. We first show that we must have
\[m_{V_I^\perp}(C_{1,\infty}) = m_{V_I^\perp}(C_{2,\infty}) = m_{V_I^\perp}(C_{\infty,\infty}) = 0\quad\quad\forall I\subseteq \bbN.\]
Combined with Lemma~\ref{lem:sometimes-excluding-4}, this justifies
restricting our attention to $Q|_{\frH_\C}$ for $\C \in
\O_{1,1}/\sim$ when calculating $\ker(Q - 4\cdot \id)$.  We will
then give that calculation, and use it to deduce the monotonicity
needed for Theorem~\ref{thm:main}.

\begin{prop}\label{prop:no-infinite-paths}
For any $I \subseteq \bbN$ we have
\[m_{V_I^\perp}(C_{1,\infty}) = m_{V_I^\perp}(C_{2,\infty}) = m_{V_I^\perp}(C_{\infty,\infty}) = 0.\]
\end{prop}

\textbf{Proof}\quad If
\[\bs{\chi} \in C_{1,\infty}\cup C_{2,\infty}\cup C_{\infty,\infty}\]
then, in particular, there is some $g \in S\cup \{e\}$ and some
singly-infinite path $P = \{g_1,g_2,\ldots\}\subseteq
\bs{\chi}^{-1}\{0\}$ starting from $g_1 \in \partial\{g\}$, and such
that for any $h \in B(P,10)\setminus P$ whose connection to $P$ does
not pass through $g$ we have $\bs{\chi}(h) = 1$.  Now given any $g_0
\in \bfF_2$ and $g_1 \in g_0 S$, let $K \subset \bfF_2$ be the
quadrant of points $h$ that are not disconnected from $g_1$ by
$g_0$. Since $B(\{e\},2)$ is finite, it will suffice to prove that
for any fixed such $g_0$ and $g_1$ we have
\begin{multline*}
m_{V^\perp}\big\{\bs{\chi}:\ \bs{\chi}^{-1}\{0\}\ \hbox{connects}\
g_1\ \hbox{to}\ \infty\ \hbox{inside}\ K\ \hbox{along some path}\ P\\
\hbox{and}\ \bs{\chi}|_{(B(P,10)\cap K)\setminus P} \equiv 1\big\} =
0.
\end{multline*}

This, in turn, will follow if we show that $m_{V^\perp}(D_N)\to 0$
as $N\to\infty$ where
\begin{multline*}
D_N:=\big\{\bs{\chi}:\ \bs{\chi}^{-1}\{0\}\ \hbox{connects}\
g_1\ \hbox{to}\ \partial B(g_1,N)\cap K\ \hbox{along some path}\ P\\
\hbox{and}\ \bs{\chi}|_{(B(P,10)\cap B(g_1,N)\cap K)\setminus P}
\equiv 1\big\}.
\end{multline*}
Now for each path $P$ that connects $g_1$ to $\partial B(g_1,N)$
inside $K$ we let \[D_{N,P}:=\big\{\bs{\chi}:\ P \subseteq
\bs{\chi}^{-1}\{0\}\ \hbox{and}\ \bs{\chi}|_{(B(P,10)\cap
B(g_1,N)\cap K)\setminus P} \equiv 1\big\},\] and now we have $D_N =
\bigcup_PD_{N,P}$.  Finally, on the one hand we know that there are
at most $3^N$ such paths $P$, and on the other we know that $P/\L_I$
has size at most $|P| = N$ for any $P$ and any $I \subseteq \bbN$,
and hence Lemma~\ref{lem:measure-equipartition} and
Corollary~\ref{cor:criterion-for-nonempty-cylinder} give
\[m_{V^\perp}(D_{N,P}) \leq \frac{1}{2^{|B(P,10)\cap B(g_1,N)\cap K| - |P/\L_I|}} \leq \frac{2^N}{2^{(2\cdot 3^9)(N-10)}}.\]
Combining these estimates gives
\[m_{V^\perp}(D_N) \leq 3^N\cdot 2^N\cdot 2^{-(2\cdot 3^9)(N-10)} \to 0\quad\quad\hbox{as}\ N\to\infty,\]
as required. \qed

\begin{cor}\label{cor:more-formulas}
For any finite path $P \subset \bfF_2$ we have
\[|\{\phi\in \bbZ_2^{B(P,10)}:\ C(\phi)\cap V_I^\perp \neq \emptyset\}| = 2^{|B(P,10)| - |P/\L_I|},\]
and so
\[\dim_{L\G_I}\ker(Q_I - 4) =
\sum_{i\geq 1}2^{-(|B(P_i,10)| - |P_i| + |P_i/\L_I|)}.\]
\end{cor}

\textbf{Proof}\quad From
Corollary~\ref{cor:criterion-for-nonempty-cylinder} and the standard
relation $|F|\cdot |F^\perp| = 2^N$ for subgroups $F \leq
\bbZ_2^N$ we can identify
\[\{\phi\in \bbZ_2^{B(P_i,10)}:\ C(\phi)\cap V_I^\perp \neq \emptyset\}\]
with the subgroup of those $\phi \in \bbZ_2^{B(P_i,10)}$ that
annihilate all the vectors of the form
\[\sum_{i=-10}^{10}\delta_{gs_1^i} - \sum_{i=-10}^{10}\delta_{ghs_1^i}\]
such that $h \in \L_I$ and $g$ and $hg$ both lie in $P_i$.  Clearly
each cell $C \in P_i/\L_I$ gives rise to a subspace of
$\bbZ_2^{B(P_i,10)}$ of dimension $|C| - 1$ spanned by these
differences with $g,hg \in C$, and so the total dimension of the
resulting subspace is
\[\sum_{C \in P_i/\L_I}(|C| - 1) = |P_i| - |P_i/\L_I|.\]

This gives the dimension of
\[\{\phi\in \bbZ_2^{B(P_i,10)}:\ C(\phi)\cap V_I^\perp \neq \emptyset\}\]
as $|B(P_i,10)| - |P_i| + |P_i/\L_I|$, and so both the desired
conclusions now follow from Lemma~\ref{lem:dim-formula}. \qed

\textbf{Proof of Theorem~\ref{thm:main}}\quad We will show that the
conclusion holds for the parameterized family of subgroups $V_I$ and
the operators $Q_I - 4$ in place of $Q_I$ provided the sequence of
lengths $l(n)$, $n\geq 1$, appearing in the definition of $t_n$
grows sufficiently fast.

Letting
\[\varphi(I) := \dim_{L\G_I}\ker(Q_I - 4)\quad\quad\hbox{for}\ I \subseteq \bbN,\]
we must prove that
\[I <_{\rm{lex}} J \quad\quad\Rightarrow\quad\quad \varphi(I) < \varphi(J)\]
provided that for each $n$ the values $l(n')$ for $n' > n$ are
sufficiently large relative to $l(1)$, $l(2)$, \ldots, $l(n)$.  More
concretely, we will prove that if $l(1) < l(2) < \ldots < l(n-1)$
and some auxiliary $L(n-1) > l(n-1)$ are such that the above
implication holds whenever $I\cap [1,n-1] <_{\rm{lex}} J\cap
[1,n-1]$ and for \emph{any} tail sequence $L(n-1)\leq l(n) < l(n+1)
< \ldots$, then we can pick particular values of $l(n)$ and $L(n)>
l(n)$ such that this same property holds with $n$ in place of $n-1$.
From this a simple recursion completes the proof.

Thus, suppose that $n\in \bbN$ is minimal such that $n\in J$ but
$n\not\in I$ (so by the definition of the lexicographic ordering we
must have $I \cap [1,n-1] = J\cap [1,n-1]$). By
Corollary~\ref{cor:more-formulas} we have
\[\varphi(J) - \varphi(I)
= \sum_{i\geq 1}\big(2^{-(|B(P_i,10)| - |P_i| + |P_i/\L_J|)} -
2^{-(|B(P_i,10)| - |P_i| + |P_i/\L_I|)}\big).\] Clearly (from the
freeness of the $t_n$s) there will be some paths $P_i$ in the above
list for which $P_i/\L_{J\cap [1,n]}$ is a nontrivial coarsening of
$P_i/\L_{J\cap [1,n-1]} = P_i/\L_{I\cap [1,n]}$ (that is, the left
cosets of the larger subgroup $\L_{J\cap [1,n]}$ intersect $P_i$ in
fewer, larger patches than those of the smaller subgroup $\L_{J\cap
[1,n-1]}$), and so the expression
\[2^{-(|B(P_i,10)| - |P_i| + |P_i/\L_{J\cap [1,n]}|)} -
2^{-(|B(P_i,10)| - |P_i| + |P_i/\L_{I\cap [1,n]}|)}\] in the above
sum will be strictly positive for each of these $i$.  Let $E \subset
\bbN$ be a finite subset of $i \in \bbN$ for which this is so, and
such that
\[\sum_{i\in E}2^{-(|B(P_i,10)| - |P_i| + |P_i/\L_{J\cap [1,n]}|)} -
2^{-(|B(P_i,10)| - |P_i| + |P_i/\L_{I\cap [1,n]}|)} =: \eta > 0.\]
Let $i_0 := \max\,E$.  Next we observe, using only very crude
estimates at every step, that for any $I$ and $J$ and any $L \geq 1$
we have
\begin{eqnarray*}
&&\sum_{i\geq 1,\,|P_i| \geq L}\big(2^{-(|B(P_i,10)| - |P_i| +
|P_i/\L_J|)} - 2^{-(|B(P_i,10)| - |P_i| + |P_i/\L_I|)}\big)\\ &&\leq
2\cdot \sum_{i\geq 1,\,|P_i|\geq L}2^{-(|B(P_i,10)| - |P_i|)} =
2\cdot \sum_{i\geq
1,\,|P_i|\geq L}2^{|P_i|}2^{-|B(P_i,10)|}\\
&&\leq 2\cdot
\sum_{\scriptsize{\begin{array}{c}\rm{all}\,\rm{paths}\,P\,\rm{in}\,\Cay(\bfF_2,S)\\
\rm{with}\,e\in B(P,1)\,\rm{and}\,|P| \geq L\end{array}}}2^{|P|}2^{-|B(P,10)|}\\
&&\leq 2\cdot
\sum_{\ell \geq L}\sum_{\scriptsize{\begin{array}{c}\rm{all}\,\rm{paths}\,P\,\rm{in}\,\Cay(\bfF_2,S)\\
\rm{with}\,e\in B(P,1)\,\rm{and}\,|P| =
\ell\end{array}}}2^{\ell}2^{-(2 + 2\cdot 3 + \ldots + 2\cdot 3^9)\ell}\\
&&\leq 2\cdot \sum_{\ell \geq L}(3\ell + 2)\cdot 3^\ell\cdot
2^\ell\cdot 2^{-(2 + 2\cdot 3 + \ldots + 2\cdot
3^9)\ell}\\
&&< \infty.
\end{eqnarray*}

Since for any finite $L$ there can be only finitely many paths among
$P_1$, $P_2$, \ldots of length $< L$, it follows that
\[\sum_{i\geq i_1 + 1}\big(2^{-(|B(P_i,10)| - |P_i| + |P_i/\L_J|)} -
2^{-(|B(P_i,10)| - |P_i| + |P_i/\L_I|)}\big)\] tends to $0$ as $i_1
\to\infty$ uniformly fast in $I$ and $J$, and so we may pick $i_1 >
i_0$ such that
\[\sum_{i\geq i_1 + 1}\big(2^{-(|B(P_i,10)| - |P_i| + |P_i/\L_J|)} -
2^{-(|B(P_i,10)| - |P_i| + |P_i/\L_I|)}\big) < \eta/2\] irrespective
of the choice of $l(n+1)$, $l(n+2)$, \ldots.  Therefore if we simply
insist that these lengths $l(n')$ for $n' > n$ should be so large
that
\[P_i/\L_I = P_i/\L_{I\cap [1,n]}\quad\quad\forall i \leq i_1,\]
we deduce that
\begin{eqnarray*}
\varphi(J) - \varphi(I) &=& \sum_{i= 1}^{i_1}\big(2^{-(|B(P_i,10)| -
|P_i| + |P_i/\L_{J\cap [1,n]}|)} - 2^{-(|B(P_i,10)| - |P_i| +
|P_i/\L_{I\cap [1,n]}|)}\big)\\
&&+ \sum_{i\geq i_1 + 1}\big(2^{-(|B(P_i,10)| - |P_i| + |P_i/\L_J|)}
-
2^{-(|B(P_i,10)| - |P_i| + |P_i/\L_I|)}\big)\\
&\geq& \sum_{i\in E}\big(2^{-(|B(P_i,10)| - |P_i| + |P_i/\L_{J\cap
[1,n]}|)} - 2^{-(|B(P_i,10)| - |P_i| + |P_i/\L_{I\cap
[1,n]}|)}\big)\\
&& - \Big|\sum_{i\geq i_1 + 1}\big(2^{-(|B(P_i,10)| - |P_i| +
|P_i/\L_J|)} - 2^{-(|B(P_i,10)| - |P_i| + |P_i/\L_I|)}\big)\Big|\\
&\geq& \eta - \eta/2= \eta/2 > 0,
\end{eqnarray*}
as required. \qed

\textbf{Remarks}\quad\textbf{1.}\quad We have presented the proof
above so as to emphasize the flexibility in the choice of the
sequence $(l(n))_{n\geq 1}$, but easy estimates show, for
example, that any doubly exponential sequence such as $l(n) =
2^{2^n}$ will do.

\quad\textbf{2.}\quad Similar arguments also prove the continuity of
the map $\phi$ for the product topology on $\P(\bbN)$, but we do not
need this. \fin

\section{Closing remarks}

Since a version of the current paper first appeared, its methods have been enhanced by Pichot, Schick and \.{Z}uk~\cite{PicSchZuk10} and Grabowski~\cite{Gra10} to obtain several further results.  On the one hand, they both show that examples answering Atiyah's question may be found among amenable groups, and with kernel dimensions equal to any chosen element of $[0,1]$.  On the other, they also both show that the kernel dimensions arising from finitely \emph{presented} groups contain all positive reals with computable binary expansions, so these also contain many non-rational examples.  The consideration of finitely presented groups has some
geometric interest since these are precisely the groups $\G$ for
which proper cocompact free $\G$-manifolds can be constructed as the
universal covers $\widetilde{M}$ of compact manifolds $M$ having
$\pi_1(M) = \G$.

Most of Pichot, Schick and \.{Z}uk's work is in refining the construction above.  On the other hand, Grabowski constructs examples by showing how to encode a program for a Turing machine into a group and group ring element.  Also, Lehner and Wagner have given even more explicit examples of non-rational algebraic von Neumann dimensions in~\cite{LehWag10}, using a construction closer to the previous paper~\cite{DicSch02}.

However, these methods seem too crude to touch what
may be the most interesting special case of Atiyah's question: that
for torsion-free groups.  This asks for a torsion-free group and a
rational group ring element $Q \in \bbQ\G$ that has any nontrivial
eigenspaces at all. It is known that this is impossible for
large classes of torsion-free groups: see, for example,~\cite{DodLinMatSchYat03},~\cite{Sch00},~\cite{Lin93} and~\cite{LinSch07}. Surveys of these results are given in~\cite{Luc02} and Reich's thesis~\cite{Rei99}.  A negative answer in general would have such striking consequences as Kaplansky's
conjecture that the group ring has no nontrivial zero-divisors.  However, I do not think the methods of the present paper offer any new hope of
constructing a positive example.

\bibliographystyle{abbrv}
\bibliography{bibfile}

\begin{thebibliography}{10}

\bibitem{Ati76}
M.~F. Atiyah.
\newblock Elliptic operators, discrete groups and von {N}eumann algebras.
\newblock In {\em Colloque ``{A}nalyse et {T}opologie'' en l'{H}onneur de
  {H}enri {C}artan ({O}rsay, 1974)}, pages 43--72. Ast\'erisque, No. 32--33.
  Soc. Math. France, Paris, 1976.

\bibitem{Aus--Atiyahexample}
T.~Austin.
\newblock Rational group ring elements with kernels having irrational
  dimension.
\newblock Preprint, available online at \verb|arXiv.org|: 0909:2360.

\bibitem{Aus--amenablepoorcompression}
T.~Austin.
\newblock Amenable groups with very poor compression into {L}ebesgue spaces.
\newblock {\em Duke Mathematical Journal}, 159(2):187--222, 2011.

\bibitem{Bol98}
B.~Bollob\'as.
\newblock {\em Modern Graph Theory}.
\newblock Springer, Berlin, 1998.

\bibitem{DenSch07}
C.~Deninger and K.~Schmidt.
\newblock Expansive algebraic actions of discrete residually finite amenable
  groups and their entropy.
\newblock {\em Ergodic Theory Dynam. Systems}, 27(3):769--786, 2007.

\bibitem{DicSch02}
W.~Dicks and T.~Schick.
\newblock The spectral measure of certain elements of the complex group ring of
  a wreath product.
\newblock {\em Geom. Dedicata}, 93:121--137, 2002.

\bibitem{DodLinMatSchYat03}
J.~Dodziuk, P.~Linnell, V.~Mathai, T.~Schick, and S.~Yates.
\newblock Approximating {$L^2$}-invariants and the {A}tiyah conjecture.
\newblock {\em Comm. Pure Appl. Math.}, 56(7):839--873, 2003.
\newblock Dedicated to the memory of J{\"u}rgen K. Moser.

\bibitem{Gra10}
L.~Grabowski.
\newblock On {T}uring dynamical systems and the {A}tiyah problem.
\newblock Preprint, available online at \verb|arXiv.org|: 1004.2030.

\bibitem{GriLinSchZuk00}
R.~I. Grigorchuk, P.~Linnell, T.~Schick, and A.~{\.Z}uk.
\newblock On a question of {A}tiyah.
\newblock {\em C. R. Acad. Sci. Paris S\'er. I Math.}, 331(9):663--668, 2000.

\bibitem{GriZuk01}
R.~I. Grigorchuk and A.~{\.Z}uk.
\newblock The lamplighter group as a group generated by a 2-state automaton,
  and its spectrum.
\newblock {\em Geom. Dedicata}, 87(1-3):209--244, 2001.

\bibitem{LehWag10}
F.~Lehner and S.~Wagner.
\newblock Free lamplighter groups and a question of {A}tiyah.
\newblock To appear, \emph{Amer. J. Math.}

\bibitem{Lin93}
P.~Linnell.
\newblock Division rings and group von {N}eumann algebras.
\newblock {\em Forum Math.}, 5(6):561--576, 1993.

\bibitem{LinSch07}
P.~Linnell and T.~Schick.
\newblock Finite group extensions and the {A}tiyah conjecture.
\newblock {\em J. Amer. Math. Soc.}, 20(4):1003--1051, 2007.

\bibitem{Luc02}
W.~L{\"u}ck.
\newblock {\em {$L\sp 2$}-invariants: theory and applications to geometry and
  {$K$}-theory}, volume~44 of {\em Ergebnisse der Mathematik und ihrer
  Grenzgebiete. 3. Folge. A Series of Modern Surveys in Mathematics [Results in
  Mathematics and Related Areas. 3rd Series. A Series of Modern Surveys in
  Mathematics]}.
\newblock Springer-Verlag, Berlin, 2002.

\bibitem{NibRol93}
G.~A. Niblo and M.~A. Roller, editors.
\newblock {\em Geometric group theory. {V}ol. 2}, volume 182 of {\em London
  Mathematical Society Lecture Note Series}, Cambridge, 1993. Cambridge
  University Press.

\bibitem{PicSchZuk10}
M.~Pichot, T.~Schick, and A.~\.{Z}uk.
\newblock Closed manifolds with transcendental ${L}^2$-{B}etti numbers.
\newblock Preprint, available online at \verb|arXiv.org|: 1005.1147.

\bibitem{Rei99}
H.~Reich.
\newblock {\em Group von Neumann algebras and related algebras}.
\newblock PhD thesis, Universit\"at G\"ottingen, 1999.
\newblock Available online from
  \href{http://www.math.uni-muenster.de/u/lueck/publ/reich/}{\texttt{http://www.math.uni-muenster.de/u/lueck/publ/reich/}}.

\bibitem{Sch00}
T.~Schick.
\newblock Integrality of {$L^2$}-{B}etti numbers.
\newblock {\em Math. Ann.}, 317(4):727--750, 2000.

\bibitem{Sch95}
K.~Schmidt.
\newblock {\em Dynamical systems of algebraic origin}, volume 128 of {\em
  Progress in Mathematics}.
\newblock Birkh\"auser Verlag, Basel, 1995.

\end{thebibliography}

\end{document}